\documentclass[11pt,reqno,tbtags]{amsart}
\usepackage{amssymb}
\usepackage[colorlinks=true, pdfstartview=FitV, linkcolor=blue,
  citecolor=blue, urlcolor=blue,pagebackref=false]{hyperref}
\usepackage[numeric,initials,nobysame]{amsrefs}
\usepackage{amsmath,amssymb,amsfonts,dsfont}
\usepackage{enumerate,graphicx}
\usepackage{algorithm}

\numberwithin{equation}{section}

\newtheorem{maintheorem}{Theorem}

\newtheorem{maincoro}{Corollary}

\newtheorem{theorem}{Theorem}[section]
\newtheorem*{theorem*}{Theorem}
\newtheorem{lemma}[theorem]{Lemma}

\newtheorem{proposition}[theorem]{Proposition}

\newtheorem{definition}[theorem]{Definition}

\theoremstyle{definition}{

}

\theoremstyle{remark}{

\newtheorem*{remark*}{Remark}

}

\newcommand{\R}{\mathbb R}
\newcommand{\N}{\mathbb N}
\newcommand{\Z}{\mathbb Z}

\newcommand{\E}{\mathbb{E}}
\renewcommand{\P}{\mathbb{P}}
\DeclareMathOperator{\var}{Var}  

\renewcommand{\epsilon}{\varepsilon}

\newcommand{\cG}{{\mathcal{G}}}
\newcommand{\cB}{{\mathcal{B}}}
\newcommand{\GPC}{\cG_{\textsc{pc}}}
\newcommand{\cF}{{\mathcal{F}}}
\newcommand{\given}{\, \big| \,}

\newcommand{\one}{\boldsymbol{1}}

\newcommand{\deq}{\stackrel{\scriptscriptstyle\triangle}{=}}
\newcommand{\K}{\mathcal{K}}
\newcommand{\GC}{{\mathcal{C}_1}} 
\newcommand{\tGC}{{\tilde{\mathcal{C}}_1}}
\newcommand{\TC}[1][\mathcal{C}_1]{#1^{(2)}} 
\newcommand{\hGC}{{\widehat{\mathcal{C}}_1}}
\newcommand{\Po}{\operatorname{Po}}
\newcommand{\Bin}{\operatorname{Bin}}
\newcommand{\Geom}{\operatorname{Geom}}
\DeclareMathOperator{\diam}{diam}
\DeclareMathOperator{\dist}{dist}

\begin{document}

\title[Anatomy of a young giant]{Anatomy of a young giant\\ component in the random graph}

\author{Jian Ding, \thinspace Jeong Han Kim, \thinspace Eyal Lubetzky and Yuval Peres}

\address{Jian Ding\hfill\break
Department of Statistics\\
UC Berkeley\\
Berkeley, CA 94720, USA.}
\email{jding@stat.berkeley.edu}
\urladdr{}

\address{Jeong Han Kim\hfill\break
Department of Mathematics, Yonsei University, Seoul 120-749 Korea, and\hfill\break
National Institute for Mathematical Sciences, Daejeon 305-340, Korea.}
\email{jehkim@yonsei.ac.kr}
\urladdr{}
\thanks{Research of J.H.\ Kim was supported by a Basic Science Research Program through the National Research Foundation of Korea (NRF) funded by the Ministry of Education, Science and Technology (CRI, No. 2008-0054850).}

\address{Eyal Lubetzky\hfill\break
Microsoft Research\\
One Microsoft Way\\
Redmond, WA 98052-6399, USA.}
\email{eyal@microsoft.com}
\urladdr{}

\address{Yuval Peres\hfill\break
Microsoft Research\\
One Microsoft Way\\
Redmond, WA 98052-6399, USA.}
\email{peres@microsoft.com}
\urladdr{}

\begin{abstract}
We provide a complete description of the giant component of the Erd\H{o}s-R\'enyi random graph $\cG(n,p)$ as soon as it emerges from the scaling window, i.e., for $p = (1+\epsilon)/n$ where $\epsilon^3 n \to \infty$ and $\epsilon=o(1)$.

Our description is particularly simple for $\epsilon = o(n^{-1/4})$, where the giant component $\GC$ is contiguous with the following model (i.e., every graph property that holds with high probability for this model also holds w.h.p.\ for $\GC$). Let $Z$ be normal with mean $\frac23 \epsilon^3 n$ and variance $\epsilon^3 n$, and let $\K$ be a random $3$-regular graph on $2\lfloor Z\rfloor$ vertices. Replace each edge of $\K$ by a path, where the path lengths are i.i.d.\ geometric with mean $1/\epsilon$. Finally, attach an independent Poisson($1-\epsilon$)-Galton-Watson tree to each vertex.

A similar picture is obtained for larger $\epsilon=o(1)$, in which case the random 3-regular graph is replaced by a random graph with $N_k$ vertices of degree $k$ for $k\geq 3$, where $N_k$ has mean and variance of order $\epsilon^k n$.

This description enables us to determine fundamental characteristics of the supercritical random graph. Namely,
 we can infer the asymptotics of the diameter of the giant component for any rate of decay of $\epsilon$, as well as the mixing time of the random walk on $\GC$.
\end{abstract}

\maketitle

\vspace{-1cm}

\section{Introduction}

The Erd\H{o}s and R\'{e}nyi random graph $\cG(n,p)$ has been studied extensively since its introduction in 1959 \cite{ER59}.
Much of the analysis of this fundamental random graph model has focused on its behavior near the critical point $p=1/n$. Nevertheless, a few key features, such as the diameter and the mixing time of the random walk on the largest component, have remained unknown in a regime just beyond criticality.

In their seminal papers from the 1960's, Erd\H{o}s and R\'enyi established a phenomenon known as \emph{the double jump}.
For $p = c/n$ where $c <1$ is fixed, the largest component $\GC$ has size $O(\log n)$ with high probability (w.h.p.).
When $c>1$, the size of $\GC$ is linear in $n$, and at the critical $c=1$ it has order $n^{2/3}$ (this latter fact
was fully established much later by Bollob\'as \cite{Bollobas84} and {\L}uczak \cite{Luczak90}).
As discovered in \cite{Bollobas84}, the critical behavior extends throughout the \emph{critical window}, the regime
where $p = (1\pm\epsilon)/n$ for $\epsilon =O(n^{-1/3})$.

Up to the critical point, the structure of $\GC$ is relatively well understood. For instance, in the fully subcritical regime ($p=(1-\epsilon)/n$ for $\epsilon > 0$ fixed), $\GC$ is a tree of known (logarithmic) size and diameter. In the critical window ($\epsilon = O(n^{-1/3}$) the distribution of $|\GC|$ was determined in \cites{LPW,Aldous2}, and the diameter was found in \cite{NP}.
See \cites{Bollobas2,JLR} for further information.

In the supercritical regime ($p=(1+\epsilon)/n$ with $\epsilon^3 n \to\infty$), a variety of methods can determine key features of $\GC$ up to some continuous functions of $\epsilon$. While these functions remain bounded in the fully supercritical case ($\epsilon > 0$ fixed), the situation becomes much more delicate as $\epsilon$ approaches the critical window.

For example, one can deduce that the diameter of the fully supercritical $\GC$ has order $\log n$ merely by analyzing certain (weak) expansion properties of its 2-core (formally defined in Section~\ref{sec-prelims}). More precise results on the diameter were obtained in \cites{RW,LS}, but they still do not give the asymptotic diameter in the whole supercritical regime.

In the fully supercritical case, it is known that the giant component consists of an expander, ``decorated'' using paths and trees of at most logarithmic size (see \cite{BKW} for a concrete example of such a statement, used there to obtain the order of the mixing time on the fully supercritical $\GC$).
However, the existing decompositions of the giant component are not precise enough to handle the case where $\epsilon \to 0$ (e.g., in \cite{RW} Riordan and Wormald point out that this is the most difficult regime for determining the diameter).

In this work, we obtain a complete characterization of the supercritical giant component.
Rather than merely describing its properties, we present a simple construction whose distribution is contiguous with that of $\GC$.
This construction is particularly elegant when the giant component is ``young'', namely when $\epsilon = o(n^{-1/4})$.
Since this is the hardest regime for alternative approaches, we start by describing this special case.

Let $\mathcal{N}(\mu,\sigma^2)$ denote the normal distribution with mean $\mu$ and variance~$\sigma^2$,
and let $\Geom(\epsilon)$ denote the geometric distribution with mean $1/\epsilon$.
\begin{maintheorem}\label{mainthm-struct}
Let $\GC$ be the largest component of the random graph $\cG(n,p)$ for $p = \frac{1 + \epsilon}{n}$, where $\epsilon^3 n\to \infty$ and $\epsilon = o(n^{-1/4})$. Then $\GC$ is contiguous to the model $\tGC$, constructed in 3 steps as follows:
\begin{enumerate}[1.\!]
  \item\label{item-struct-base} Let $Z \sim \mathcal{N}\left(\tfrac23\epsilon^3 n, \epsilon^3 n\right)$, and 
      select a random 3-regular multigraph $\K$ on $N = 2\lfloor Z \rfloor$ vertices.

  \item\label{item-struct-edges} Replace each edge of $\K$ by a path, where the path lengths are i.i.d.\ $\Geom(\epsilon)$.
  \item\label{item-struct-bushes} Attach an independent $\mathrm{Poisson}(1-\epsilon)$-Galton-Watson tree to each vertex.
\end{enumerate}
That is, $\P(\tGC \in \mathcal{A}) \to 0$ implies $\P(\GC \in \mathcal{A}) \to 0$
for any set of graphs $\mathcal{A}$.
\end{maintheorem}
In the above, a $\mathrm{Poisson}(\mu)$-Galton-Watson tree is the family tree of a Galton-Watson branching process with offspring distribution $\mathrm{Poisson}(\mu)$.

Two well-known objects relevant to the study of the giant component are its 2-core $\TC$ and its kernel $\K$. The 2-core of a graph is its maximum subgraph where all degrees are at least $2$. The kernel is obtained from the 2-core by replacing every maximal 2-path by an edge (where a 2-path is a path where all internal vertices have degree 2).
Note that our description of $\tGC$ constructs the kernel in Step~\ref{item-struct-base}, the 2-core in Step~\ref{item-struct-edges} and the entire component $\tGC$ in Step~\ref{item-struct-bushes}.

The above theorem not only states that the kernel of $\GC$ in this regime is an expander, but it is in fact contiguous to a random 3-regular graph, an object whose expansion properties are well understood (cf., e.g., \cite{HLW}).
Furthermore, the 2-core is obtained from the kernel by a simple operation (``stretching'' the edges into paths of lengths i.i.d.\ geometric with mean $1/\epsilon$). This allows us to pinpoint the expansion properties of the 2-core and their dependence on $\epsilon$ as it tends to $0$.

A few known (yet nontrivial) properties of the 2-core of $\GC$ can be immediately read off from Theorem~\ref{mainthm-struct}. For instance, w.h.p.\ the 2-core contains $(2+o(1))\epsilon^2 n$ vertices while the kernel has $(\tfrac43+o(1))\epsilon^3 n$ vertices (see \cites{Luczak91,PW}).
As there are w.h.p.\ $(2+o(1))\epsilon^3 n$ edges in the kernel, a simple estimate of the maximum of i.i.d.\ geometric variables gives the following corollary.

\begin{maincoro}\label{maincoro-2-core-2-path} Let $\TC$ be the $2$-core of the largest component of $\cG(n,p)$ for $p = \frac{1 + \epsilon}{n}$, where $\epsilon^3 n\to \infty$ and $\epsilon = o(n^{-1/4})$. The maximal 2-path in $\TC$ has length
$(1/\epsilon)\log(\epsilon^3 n) + O_\mathrm{P}(1/\epsilon)$.
\end{maincoro}
Similarly, since a random $3$-regular graph is Hamiltonian w.h.p.\ (see \cite{RoW}), we immediately deduce that in the above regime $\GC$ contains a simple cycle of length $(\frac43+o(1))\epsilon^2 n$. This matches the lower bound of {\L}uczak \cite{Luczak91} on the circumference of the supercritical random graph.

Moreover, Theorem~\ref{mainthm-struct} enables us to interpret distances in the 2-core as passage times in first-passage percolation (for further information on this thoroughly studied topic, see, e.g., \cite{Kesten}). As we state in Theorem~\ref{mainthm-diam-gc} below, this connection (used in a companion paper \cite{DKLP}) gives the asymptotic behavior of the diameter throughout the regime $\epsilon^3 n \to \infty$ and $\epsilon=o(1)$.

\subsection{Main results}

We now state the extension of Theorem~\ref{mainthm-struct} to all $\epsilon=o(1)$ outside the critical window.

\begin{maintheorem}\label{mainthm-struct-gen}
Let $\GC$ be the largest component of $\cG(n,p)$ for $p = \frac{1 + \epsilon}{n}$, where $\epsilon^3 n\to \infty$ and $\epsilon\to 0$. Let $\mu<1$ denote the conjugate of $1+\epsilon$, that is,
$\mu\mathrm{e}^{-\mu} = (1+\epsilon) \mathrm{e}^{-(1+\epsilon)}$. Then $\GC$ is contiguous to the following model $\tGC$:
\begin{enumerate}[\indent 1.]
  \item\label{item-struct-gen-degrees} Let $\Lambda\sim \mathcal{N}\left(1+\epsilon - \mu, \frac1{\epsilon n}\right)$ and assign
  i.i.d.\  variables $D_u \sim \mathrm{Poisson}(\Lambda)$ ($u \in [n]$) to the vertices, conditioned that $\sum D_u \one_{D_u\geq 3}$ is even.
  Let
  \[\mbox{$N_k = \#\{u : D_u = k\}$ \quad and \quad $N= \sum_{k\geq 3}N_k$}\,.\]  
  Select a random multigraph $\K$ on $N$ vertices, uniformly among all multigraphs with $N_k$ vertices of degree $k$ for $k\geq 3$.
  \item\label{item-struct-gen-edges} Replace the edges of $\K$ by paths of lengths i.i.d.\ $\Geom(1-\mu)$. 
  \item\label{item-struct-gen-bushes} Attach an independent $\mathrm{Poisson}(\mu)$-Galton-Watson tree to each vertex.
\end{enumerate}
That is, $\P(\tGC \in \mathcal{A}) \to 0$ implies $\P(\GC \in \mathcal{A}) \to 0$
for any set of graphs $\mathcal{A}$.
\end{maintheorem}
We note that conditioning that the sum of degrees is even can easily be realized by rejection sampling.
The differences between the two theorems are the approximation of $1-\mu \approx \epsilon$ in Steps~2,3, and a richer degree distribution of the random graph $\K$ in Step~1. 

Further note that it was shown by {\L}uczak \cite{Luczak91} that the kernel $\K$ in the above regime is a random multigraph on a certain degree sequence, which is cubic except for a negligible number of vertices. However, in that description the vertex degrees and lengths of the 2-paths subdividing the kernel edges are all dependent, whereas in our contiguous model these are i.i.d.\ Poisson (Step~1) and i.i.d.\ Geometric (Step~2) respectively.


Combining Theorem~\ref{mainthm-struct-gen} with some known results on first-passage percolation from \cite{BHvdH} gives an immediate corollary on the typical distances between vertices of degree at least 3 in the 2-core.
\begin{maincoro}\label{maincoro-kernel-typical-dist} Let $\TC$ be the $2$-core of the largest component of $\cG(n,p)$ for $p = \frac{1 + \epsilon}{n}$, where $\epsilon^3 n\to \infty$ and $\epsilon = o(1)$. Let $u,v$ be two vertices of degree at least $3$ in $\TC$, chosen u.a.r.\ among all such vertices. The distance between $u,v$ is w.h.p.\ $(1/\epsilon + O(1))\log(\epsilon^3 n)$.
\end{maincoro}

However, maximal distances in the 2-core can differ from typical distances; compare the above result to \eqref{eq-diam-kernel} in the next theorem, which we prove in a companion paper.

\begin{maintheorem}[\cite{DKLP}] \label{mainthm-diam-gc}
Consider the random graph $\cG(n,p)$ for $p = \frac{1 + \epsilon}{n}$, where $\epsilon^3 n\to \infty$ and $\epsilon = o(1)$.
Let $\GC$ be the largest component $G$, let $\TC$ be its $2$-core and let $\K$ denote its kernel.
Then w.h.p.,
\begin{align}
\diam(\GC) &= \big(3+o(1)\big)(1/\epsilon)\log(\epsilon^3 n)\,,\label{eq-diam-gc}\\
\diam(\TC) &= \big(2+o(1)\big)(1/\epsilon)\log(\epsilon^3 n)\,,\label{eq-diam-2-core}\\
  \max_{u,v\in \K} \dist_{\TC}(u,v) &= \big(\tfrac53+o(1)\big)(1/\epsilon)\log (\epsilon^3 n)\label{eq-diam-kernel}\,.
\end{align}
\end{maintheorem}

To prove the above theorem, we need to go beyond typical distances and obtain new large deviation estimates for the relevant parameters (see \cite{DKLP} for further details). The result \eqref{eq-diam-gc} on the diameter of the giant component concludes a long list of studies of this parameter in the supercritical random graph (e.g., \cites{CL,FeRa,LS,RW}). First results for the challenging regime where $\epsilon=o(1)$ appeared only recently: Riordan and Wormald \cite{RW} obtained very accurate estimates of the diameter for most of this regime, but did not cover the range where the random graph emerges from the critical window (i.e., $\epsilon^3 n$ tends to $\infty$ arbitrarily slowly). {\L}uczak and Seierstad \cite{LS} then gave estimates for the diameter that do apply to the entire supercritical regime, yet their upper and lower bounds differ by a factor of $\frac{1000}7$.

Controlling typical and maximal distances between vertices in the giant component is but one of several prerequisites for estimating the mixing time of the (lazy) random walk on $\GC$. For instance, as this parameter is highly sensitive to bottlenecks in $\GC$, one also needs to fully understand the isoperimetric profile of the 2-core and the structure of the trees attached to it.

In the fully supercritical case, Fountoulakis and Reed \cite{FR} and Benjamini, Kozma and Wormald \cite{BKW} independently proved that
the mixing time on $\GC$ is of order $\log^2 n$. However, as evident from the structure description in Theorem~\ref{mainthm-struct-gen}, methods for the fully-supercritical case that depend on large sets in the 2-core having edge expansion bounded away from $0$ will break down as $\epsilon \to 0$. Within the critical window, it was shown in \cite{NP} that the mixing time on $\GC$ has order $n$.
For $\epsilon=o(1)$ outside the critical window, the problem of estimating the mixing-time on $\GC$ remained open, and furthermore, it was unclear what the answer should be, as one would expect some interpolation between $\log^2 n$ for fixed $\epsilon>0$ and
order $n$ at criticality.

The following theorem, proved in a companion paper, settles this problem by exploiting the geometric understanding of $\GC$ provided by Theorem~\ref{mainthm-struct-gen}. This completes the picture of the supercritical mixing time.

\begin{maintheorem}[\cite{DLP}]\label{mainthm-tmix}
Let $\GC$ be the largest component of $\cG(n,p)$ for $p = \frac{1 + \epsilon}{n}$, where $\epsilon^3 n\to \infty$ and $\epsilon=o(1)$.
With high probability, the mixing time of the lazy random walk on $\GC$ is of order $(1/\epsilon^3) \log^2 (\epsilon^3 n)$.
\end{maintheorem}

Indeed, the mixing time exhibits a smooth evolution from the critical regime $\epsilon = O(n^{-1/3})$ to the fully-supercritical regime of $\epsilon>0$ fixed.

\subsection{Main techniques}
A key ingredient in the proofs is the Poisson cloning model $\GPC(n,p)$, introduced in \cite{KimB} and shown to be contiguous to $\cG(n,p)$ (see Section~\ref{sec-prelims}). It thus suffices to establish the contiguity of our model $\tGC$ to the giant component of Poisson cloning, a fact we establish in several stages.

We first show the contiguity of the 2-cores in the models through a careful analysis of $\GPC(n,p)$. We then perform a series of contiguous translations of the model, in order to remove dependencies between maximal 2-paths in the 2-core, as well as incorporate the trees attached to the 2-core in $\GC$. To establish these, we use local central limit results for various parameters, including a powerful local CLT of Pittel and Wormald \cite{PW}.

\subsection{Organization}
 Section~\ref{sec-prelims} contains several preliminary facts needed for the proofs. In Section~\ref{sec:structure-2-core} we reduce the 2-core of the Poisson cloning model to an intermediate simplified model (the proof of a technical lemma on Poisson cloning used here is postponed to
Section~\ref{sec:lambda-C}). This model is subsequently reduced in Section~\ref{sec:poisson-geometric} to one that is essentially the 2-core of our model $\tGC$. The complete structure of the giant component is thereafter analyzed in Section~\ref{sec:structure-gc}, which concludes the proof of Theorem~\ref{mainthm-struct-gen}. In Section~\ref{sec:simplify-model} we prove Theorem~\ref{mainthm-struct}, addressing the special case of the early giant component.


\section{Preliminaries}\label{sec-prelims}

\subsection{Cores and kernels} The \emph{$k$-core} of a graph $G$, denoted by $G^{(k)}$, is its maximum subgraph $H \subset G$ where every vertex has degree at least $k$. It is well known (and easy to verify) that this subgraph is unique, and can be obtained by repeatedly deleting any vertex whose degree is smaller than $k$ (at an arbitrary order).

We call a path $\mathcal{P} = v_0, v_1, \ldots, v_k$ for $k > 1$ (i.e., a sequence of vertices with $v_i v_{i+1}$ an edge for each $i$) a \emph{2-path} if and only if $v_i$ has
  degree $2$ for all $i = 1, \ldots, k-1$ (while the endpoints $v_0,v_k$ may have degree larger than $2$, and possibly $v_0=v_k$).

The \emph{kernel} $\K$ of $G$ is obtained by taking its $2$-core $\TC[G]$ minus its disjoint cycles, then repeatedly contracting all 2-paths (replacing each by a single edge). Notice that, by definition, the degree of every vertex in $\K$ is at least $3$.
At certain times the notation $\ker(G)$ will be useful to denote a kernel with respect to some specific graph $G$.

\subsection{Configuration model}\label{sec-prelim-conf}
This model, introduced by Bollob\'{a}s \cite{Bollobas1}, provides a remarkable method for constructing random graphs with a given degree distribution, which is highly useful to their analysis. We describe this for the case of random $d$-regular graphs for $d$ fixed (the model is similar for other degree distributions); see \cites{Bollobas2,JLR,Wormald99} for additional information.

Associate each of the $n$ vertices with $d$ distinct points (also referred to as ``half-edges''), and consider a uniform perfect matching on these points. The random $d$-regular graph is obtained by contracting each cluster of the $d$ points corresponding to a vertex, possibly introducing multiple edges and self-loops. Clearly, on the event that the obtained graph is simple, it is uniformly distributed among all $d$-regular graphs, and furthermore, one can show that this event occurs with probability bounded away from $0$ (namely, with probability about $\exp(\frac{1-d^2}4)$). Hence, every event that occurs w.h.p.\ for this model, also occurs w.h.p.\ for a random $d$-regular graph.

One particularly useful property of the above model is that it allows one to construct the graph gradually, exposing the edges of the matching one by one. This way, having exposed part of the graph, the edges on the remaining unmatched points are still distributed as a uniform perfect matching.

\subsection{Poisson cloning model}\label{sec-prelim-pc}
In order to analyze the delicate structure of the near-critical giant component, we need to use Poisson cloning model $\GPC(n,p)$, which was introduced in \cite{KimB}. We incorporate a brief account on Poisson cloning model as follows, and one can see \cite{KimA} and \cite{KimB} for more.

Let $V$ be the set of $n$ vertices, and $\Po(\lambda)$ denote a Poisson random variable with mean $\lambda$.
Let $\{d_v\}_{v \in V}$ be a sequence of i.i.d.\ $\Po(\lambda)$ variables with $\lambda = (n-1)p$. Then, take $d(v)$ copies of each vertex $v \in V$ and the copies of $v$ are called \emph{clones} of $v$ or simply \emph{$v$-clones}. Define $N_\lambda \deq \sum_{v \in V} d(v)$.

If $N_\lambda$ is even, the multi-graph $\GPC(n,p)$ is obtained by generating a uniform random perfect matching of those $N_\lambda$ clones (e.g., via the configuration model, where every clone is considered to be a half-edge) and contracting clones of the same vertex. That is to say, each matching of a $v$-clone and a $w$-clone is translated into the edge $(v, w)$ with multiplicity. In the case that $v = w$, it contributes a self-loop with degree $2$. On the other hand, if $N_\lambda$ is odd, we first pick a uniform clone and translate it to a special self-loop contributing degree $1$ of the corresponding vertex. For the remaining clones, we generate a perfect matching and contract them as in the $N_\lambda$ even case.

The following theorem of \cite{KimB} states that the Poisson cloning model is \emph{contiguous} with Erd\H{o}s-R\'{e}nyi model. Hence, it suffices to study Poisson cloning model in order to establish properties of Erd\H{o}s-R\'{e}nyi model.

\begin{theorem}[\cite{KimB}*{Theorem 1.1}]\label{thm-poisson-ER}
Suppose $p = \Theta (n^{-1})$. Then there exist constants $c_1, c_2 >0$ such that for any collection $\cF$ of simple graphs, we have
$$c_1\P(\cG_{\textsc{pc}}(n,p) \in \cF) \leq \P(\cG(n,p) \in \cF) \leq c_2\big(\P(\GPC(n,p) \in \cF)\big)^{1/2} + \mathrm{e}^{-n}\big)~.$$
\end{theorem}

Note that in our regime ($p = \frac{1+\epsilon}n$ for $\epsilon=o(1)$ and $\epsilon^3 n \to \infty$) we may replace the rate $\lambda=(n-1)p$ in the Poisson-cloning model definition simply by $\lambda=np$, for convenience.

\section{The 2-core of Poisson cloning}\label{sec:structure-2-core}

By the results of \cite{KimB}, the random graph $\cG(n,p)$ in our range of parameters is contiguous to the Poisson cloning model, where every vertex gets an i.i.d.\ $\Po(np)$ number of half-edges (clones), and the final (multi)graph is obtained thereafter via the configuration model.
As opposed to $\cG(n,p)$, the Poisson cloning model features vertex degrees that are independently distributed, often contributing to an easier analysis. Nevertheless, the structure of the 2-core in this model just beyond criticality is still highly nontrivial.

The main goal in this section is to reduce the 2-core of the supercritical Poisson cloning model to the following tractable model, which is simply a random graph uniformly chosen over all graphs with a given degree sequence.

\begin{definition}[Poisson-configuration model for $n$ and $p=\frac{1+\epsilon}n$]\label{def-poisson-conf-model}\mbox{}
\begin{enumerate}[(1)]
\item Let $\Lambda\sim \mathcal{N}\left(1+\epsilon - \mu, \frac1{\epsilon n}\right)$ and assign an independent variable $D_u \sim \Po(\Lambda)$ to each vertex $u$. Let $N_k = \#\{u : D_u = k\}$ and $N = \sum_{k\geq 2}N_k$.
\item Construct a random graph on $N$ vertices, uniformly chosen over all graphs with $N_k$ degree-$k$ vertices for $k\geq 2$ (if $N$ is odd, choose a vertex $u$ with $D_u=k \geq 2$ with probability proportional to $k$, and give it $k-1$ half-edges and a self-loop).
\end{enumerate}
\end{definition}

\begin{theorem}\label{thm-Lambda-contiguity}
Let $G \sim \GPC(n,p)$ be generated by the Poisson cloning model for
$p=\frac{1+\epsilon}n$, where $\epsilon\to 0$ and $\epsilon^3 n \to \infty$.
Let $\TC[G]$ be its $2$-core, and $H$
be generated by the Poisson-configuration model corresponding to $n,p$. Then for any set of graphs $\mathcal{A}$ such that $\P(H \in
\mathcal{A}) \to 0$, we have $\P(\TC[G] \in \mathcal{A}) \to 0$.
\end{theorem}


In order to prove the above Theorem~\ref{thm-Lambda-contiguity}, in what follows we review a specific way to generate $\GPC(n,p)$, introduced in \cite{KimB}. Let $V$ be a set of $n$ vertices and let
$$\lambda \deq  n p = 1+ \epsilon~,$$
be the mean of the degree. Consider $n$ horizontal line segments ranging from $(0, j)$ to $(\lambda, j)$, for $j=1, \ldots, n$ in $\R^2$. Assign a Poisson point process with rate $1$ on each line segment independently. Each point $(x, v)$ in these processes is referred to as a $v$-clone with the assigned number $x$. The entire set of Poisson point processes is called a Poisson $\lambda$-cell.

Given the Poisson $\lambda$-cell, there are various schemes to generate a perfect matching on all points (thus yielding a random graph). One such way is the ``Cut-Off Line Algorithm'' (COLA), defined in \cite{KimA}, which is useful in finding the 2-core $\TC[G]$.
We next describe this algorithm in detail.

First define $\theta_\lambda$ to be the unique positive solution to the following equation:
\begin{equation}\label{eq-theta-lambda-def}
\theta = 1 - \mathrm{e}^{-\theta \lambda}~.
\end{equation}
It is straightforward to verify that
\begin{equation}\label{eq-theta-lambda}
\theta_\lambda = (2+ o(1)) \epsilon~.
\end{equation}
Let $\beta$ be some real, to be specified later, satisfying
\begin{equation}
  \label{eq-beta-req1}
\frac{1- \theta_\lambda}{3} \leq \beta \leq \frac{1-\theta_\lambda}{2}~.
\end{equation}
Next, construct a Poisson $\lambda$-cell as follows. The COLA procedure consists of multiple phases, formally defined in Algorithm~\ref{algorithm-cola} below. Throughout these phases, the algorithm maintains the position of a ``cut-off line'', a vertical line in $\R^2$ whose initial $x$-coordinate equals $\lambda$, and gradually moves leftwards. The $j$-th phase ($j \geq 1$) begins when the line is at $(1-\beta)^{j-1}\lambda$ and
ends once it reaches $(1-\beta)^j \lambda$.

The result of each phase is a matching on (previously unmatched) clones. In order to describe the rule of constructing this matching, we need the following definitions. At any given point, we call a vertex $v \in V$ (and its unmatched clones) \emph{light} if it has at most one unmatched clone and \emph{heavy} otherwise. Furthermore, for each $j$, we label each vertex (and its unmatched clones) at the beginning of phase $j$ as either $j$-\emph{active} or $j$-\emph{passive}, as follows. A vertex $v\in V$ (and its clones) is $j$-passive if
it has precisely $2$ clones to the left of the cut-off line, and both are unmatched. This partition of the unmatched clones into $j$-active and $j$-passive ones remains fixed throughout phase $j$.

At the beginning of the process, all the light clones are placed in a stack (whose state is maintained without being re-initialized after each phase). The order by which these clones are inserted into the stack can be arbitrary, as long as it is oblivious of the values assigned to the clones.

\begin{algorithm}
$$
\fbox{\parbox{.98\linewidth}{\mdseries
\begin{enumerate}[1.]
\item \label{item-cola-1} As long as the stack is nonempty, repeat the following:
\begin{itemize}
\item Let $(u,i)$ be the first clone in the stack.
\item Move the cut-off line leftwards until one of the following occurs:
\begin{enumerate}
  \item If the line hits $(1-\beta)^j$, the phase is conlcuded (quit).
  \item The line hits an unmatched clone $(v,j) \neq (u,i)$.
\end{enumerate}
\item Remove $(u,i)$ from the stack, as well as $(v,j)$ (if it is there).
\item Match $(u,i)$ and $(v,j)$, and re-evaluate $u$ and $v$ as light/heavy.
\item Add any clone that just became \emph{light} into the stack.
\end{itemize}
\item \label{item-cola-2} If there are \emph{active} unmatched clones:
 \begin{enumerate}
   \item\label{item-cola-2-1}
   Choose such a clone uniformly at random and put it in the stack.
   \item Return to Step~\eqref{item-cola-1}.
 \end{enumerate}
\noindent Otherwise, the algorithm is concluded (no additional phases).
\end{enumerate}
}}
$$
\caption{\textsc{Cut-Off Line Algorithm: phase $j$ description}}
\label{algorithm-cola}
\end{algorithm}
Define $\Lambda_C$ to be the $x$-coordinate of the cut-off line once Step~\ref{item-cola-2} is reached for the first time in the course of the algorithm, i.e., at the first time when there are no light clones. The next lemma states that $\Lambda_C$ is concentrated about
$\theta_\lambda \lambda$ with a standard deviation of $1/\sqrt{\theta_\lambda n}$.

Before giving the explicit statement on the concentration of $\Lambda_C$, we elaborate on its important role in understanding the structure of the 2-core of the graph. Until reaching Step~\ref{item-cola-2} for the first time, the above algorithm
repeatedly matches light clones until all of them are exhausted --- precisely as the cut-off line reaches $\Lambda_C$. As stated in Section~\ref{sec-prelims}, the $k$-core of a graph can be obtained by repeatedly removing vertices of degree at most $k-1$ (at any arbitrary order). Therefore, the $2$-core is precisely comprised of all the unmatched clones at the moment we reach $\Lambda_C$.
Crucially, continuing the algorithm will further reveal the inner structure of the $2$-core, and these further steps are equivalent to running the configuration model on the clones to the left of $\Lambda_C$.

The following theorem gives tight concentration bounds for $\Lambda_C$. Its proof relies on a delicate analysis of the above mentioned
Algorithm~\ref{algorithm-cola}, and we postpone it to Section~\ref{sec:lambda-C}.
\begin{theorem}\label{thm-Lambda-C-upper}[Upper bound on the window of $\Lambda_C$]
There exist some constant $c > 0$ so that for all
$\gamma>0$ with $\gamma = o\big(\sqrt{\theta_\lambda^3 n}\big)$, the following holds:
\begin{equation}\label{eq-Lambda-C-window}
\P\left(|\Lambda_C - \theta_\lambda \lambda| \geq \tfrac{\gamma}{\sqrt{\theta_\lambda n}} \right) \leq
\mathrm{e}^{-c \gamma^2}\,.
\end{equation}
\end{theorem}

\subsection{Size of the 2-core and its disjoint cycles}
Using the above theorem, we will now be able to characterize the structure of $\TC[G]$.
Indeed, by the discussion preceding the theorem, the $2$-core of Poisson-cloning given that $\Lambda_C=\ell$ has the same distribution as the graph generated by the Poisson-configuration model given that $\Lambda=\ell$. The above theorem implies that w.h.p.\ we only need to consider \[\ell = (1+o(1))\theta_\lambda\lambda = (2+o(1))\epsilon\,,\] and our next step is to estimate the basic properties of the 2-core (size, the number of vertices that comprise disjoint cycles) and its kernel on this event.

The next proposition thus applies not only to the Poisson-configuration model but also to the 2-core of Poisson-cloning.
The term \emph{expander} used here refers (informally) to a graph where the ratio between the boundary and volume of each set is bounded from below by some constant $c>0$ (a precise definition appears below).

\begin{proposition}\label{prop-2-core-structure}
Let $H$ be generated by the Poisson-configuration model given $\Lambda = \ell$, where $\ell=(2+o(1))\epsilon$.
Define $H'$ as the graph obtained by deleting every disjoint cycle from $H$.
Let $N_2$ be the number of vertices with degree $2$ in $H$, and $N'_2$ be the corresponding quantity for $H'$.
Then w.h.p.
\begin{align*}
N_2 &= (2+o(1)) \epsilon^2 n, & N'_2 = (1+o(1))N_2\,.
\end{align*}
In addition, w.h.p.\ the kernel $\K$ of $H$ is an expander graph with
\begin{align*}
|\K| &= \left(\tfrac{4}{3}+ o(1)\right) \epsilon^3 n~,& |E(\K)| = (2+o(1))\epsilon^3 n\,.
\end{align*}
\end{proposition}


The first step in the proof is to establish the size of the kernel $\K=\ker(H)$, as well as show that it is an expander. This latter fact is of independent interest and will have important applications, e.g., for the mixing time of the random walk on $\GC$.
In what follows, for a subset $S$ of the vertices of a graph $G$, we let
$$ d_G(S) \deq \sum_{v\in S} d_G(v)$$
denote the sum of the degrees of its vertices (also referred to as the \emph{volume} of $S$ in $G$). Further define the \emph{isoperimetric number} of a graph $G$, denoted by $i(G)$,
as
$$ i(G) \deq \min \left\{ \frac{e(S,S^c)}{d_G(S)}\;:\; S \subset V(G)\,,\, d_G(S) \leq |E(G)| \right\}~.$$
We say that $G$ is a $c$-\emph{edge-expander} for some fixed $c> 0$ iff $i(G) > c$.

\begin{lemma}\label{lem-kernel-expander}
For $\K$ the kernel of $H$ as defined in Proposition~\ref{prop-2-core-structure}, w.h.p.
\[|\K|=\left(\tfrac43+o(1)\right)\epsilon^3 n~,\quad|E(\K)| = (2+o(1))\epsilon^3 n~,\]
and $\K$ is an $\alpha$-edge-expander for some constant $\alpha>0$.
\end{lemma}
\begin{proof}
By Definition~\ref{def-poisson-conf-model}, the kernel consists of exactly those vertices $u\in V$ that have $D_u \geq 3$. Combining this with
the assumption that $\Lambda = (2 + o(1)) \epsilon$, it follows that
$|\K| \sim \Bin(n, p^+_{3}(\Lambda))$, where
\begin{equation*}
p^+_{3}(\Lambda) \deq \sum_{k\geq 3} \mathrm{e}^{-\Lambda} \frac{\Lambda^k}{k!} = \mathrm{e}^{-\Lambda}(1 + O(\Lambda))\frac{\Lambda^3}{3!} = \left(\tfrac{4}{3} + o(1)\right) \epsilon^3~.
\end{equation*}
Since $\epsilon^3 n \to \infty$, we get that w.h.p.
\begin{equation}\label{eq-kernel-size}|\K| = \left(\tfrac{4}{3}+ o(1)\right) \epsilon^3 n~.\end{equation}
Similarly, the total sum of degrees in $\K$ is simply $\sum_{u \in \K} D_u$, and therefore, $|E(\K)|$ is the sum of $n$ i.i.d.\ variables distributed as $Y \sim \Po(\Lambda)\one_{[3,\infty)}$. A similar calculation to the one above now gives
that $\E Y = (4+o(1))\epsilon^3$, and so (by CLT) w.h.p.
\begin{equation}\label{eq-kernel-edge-size}
|E(\K)| = (2+ o(1)) \epsilon^3 n~.
\end{equation}

To show that the isoperimetric number is bounded away from $0$, we apply standard techniques used to analyze the configuration model (for definitions, see Subsection~\ref{sec-prelim-conf}), while assuming \eqref{eq-kernel-size} and \eqref{eq-kernel-edge-size}.
Let $0 < \alpha < \frac14$ be specified later, and let $$D = 2|E(\K)| = (4+o(1))\epsilon^3 n$$ be the total number of points to be matched in the configuration model. We will next prove a lower bound $i(\K) > \alpha > 0$ for the case where $D$ is even, and under the relaxed condition that perhaps one of the vertices of $\K$ has degree $2$ (all others have degree $3$ or more).

To see that this gives a bound on $i(\K)$ for $D$ odd, recall that in that case precisely one of the $D$ points will have a self-loop (by Definition~\ref{def-poisson-conf-model}). Clearly, omitting one point produces a kernel as handled above (with an isoperimetric number at least $\alpha$), and reintroducing it (to a vertex with at least $2$ other points) would give $i(\K) \geq \frac23 \alpha$.

 Consider the probability that $\{e(S, S^c) \leq \alpha d_\K(S)\}$, where $S$ is a subset of $\K$ with $d_\K(S) = s$. This is precisely the probability that $k \leq \alpha s$ points out of the $s$ that comprise $S$ are matched with points in $S^c$, whereas the remaining $s-k$ points form a perfect matching.
Thus,
\begin{align*}
\P(e(S, S^c) &\leq \alpha s) = \sum_{k = 0}^{\alpha s} \P(e(S, S^c) = k)\one_{\{s \equiv k \!\!\!\!\pmod{2}\}}\\
&\leq\frac{1}{(D-1)!!}\sum_{k = 0}^{\alpha s}\binom{s}{k}(s-k-1)!! \binom{D - s}{k}k! (D - s - k - 1)!!\\
&\leq \frac{4}{D!!}\sum_{k = 0}^{\alpha s}\binom{s}{k}(s-k)!! \binom{D - s}{k}k! (D - s - k)!!~,
\end{align*}
where we used the fact that $\frac{m!!}{(m-1)!!} < 2\sqrt{m}$ for sufficiently large $m$,
and that $\sqrt{(D-s-k)(s-k)} \geq \sqrt{D-2k} \geq \sqrt{D/2}$, as $k \leq \alpha s \leq s/4$ and $s \leq D/2$.

A standard application of Stirling's formula gives $n!!  = \Theta ((n!)^{1/2} n^{1/4})$. Hence, for some constant $c_1 > 0$,
\begin{align*}
\P(e(S, S^c) \leq \alpha s) &\leq c_1 \sum_{k = 0}^{\alpha s} \frac{(D - s-k)^{1/4}}{(D!)^{1/2}D^{1/4}}\frac{s!(s-k)^{1/4}}{k! ((s-k)!)^{1/2}} \frac{(D - s)!}{((D - s - k)!)^{1/2} } \\
&= c_1 \sum_{k = 0}^{\alpha s} \bigg(\frac{(D-s-k)(s-k)}{D}\bigg)^{1/4}\bigg(\frac{\binom{s}{k}\binom{D-s}{k}}{\binom{D}{s}}\bigg)^{1/2}\\
&\leq c_1 s^{1/4} \sum_{k = 0}^{\alpha s} \bigg(\frac{\binom{s}{k}\binom{D-s}{k}}{\binom{D}{s}}\bigg)^{1/2}~.
\end{align*}
It is well known that, by Stirling's formula, $\binom{n}{k} \asymp \sqrt{\frac{n}{k(n-k)}}\exp\left[-H(\frac{k}n)n\right]$,
where $H(x)$ is the entropy function $H(x) \deq - x\log x - (1-x) \log(1-x)$.
Thus, for some constant $c_2>0$,
\begin{align*}
\P(e(S, S^c) \leq \alpha s)
&\leq c_2 \sqrt{s} \sum_{k=0}^{\alpha s} \mathrm{e}^{\frac{1}{2}\left[H\left(\frac{k}{s}\right)s + H\left(\frac{k}{D-s}\right)(D - s)- H\left(\frac{s}{D}\right)D\right]}\\
&\leq c_2 s^{3/2} \mathrm{e}^{\frac{1}{2}\left[H\left(\alpha\right)s + H\left(\frac{2\alpha s}{D}\right)D - H\left(\frac{s}{D}\right)D\right]}~,
\end{align*}
where we applied the fact that $\frac{k}{D-s} \leq \frac{\alpha s}{D/2}$.
Recalling that each vertex of $\K$ has degree at least $3$, except for possibly one vertex of degree $2$, we have
\[s = d_\K(S) \geq 3 |S|-1~,\quad D \geq 3 |\K|-1~,\] and it follows that $|S| \leq \frac{s+1}3$. Thus,
\begin{align*}
\sum_{S: d_\K(S) = s}&\P(e(S, S^c) \leq \alpha s) \leq \sum_{l \leq \frac{s+1}3}\sum_{\substack{S: d_\K(S)=s \\ |S|=l }}\P(e(S, S^c) < \alpha s)\\
& \leq \sum_{l\leq \frac{s+1}3} c_2 s^{3/2} \binom{|\K|}{l} \mathrm{e}^{\frac{1}{2}\left[H\left(\alpha\right)s + H\left(\frac{2\alpha s}{D}\right)D - H\left(\frac{s}{D}\right)D\right]}~.
\end{align*}
Since $D = (3+o(1))|\K|$, any $s \leq D/2$ satisfies
$s/3 \leq (\frac12 + o(1)) |\K|$. Therefore, another application of the above estimate of the binomial coefficient gives that
for some constant $c_3> 0$,
\begin{align*}
\sum_{S: d_\K(S) = s}\P(e(S, S^c) \leq \alpha s) &\leq c_3 s^{5/2} \mathrm{e}^{H(\frac{s}{3|\K|} +o(1))|\K|+ \frac{1}{2}\left[H(\alpha)s + H(\frac{2 \alpha s}{ D})D - H(\frac{s}{D})D\right]}\\
&= c_3 s^{5/2} \mathrm{e}^{(\frac13+o(1))H(\frac{s}{D})D + \frac{1}{2}\left[H(\alpha)\frac{s}D + H(\frac{2 \alpha s}{ D}) - H(\frac{s}{D})\right]D}~.
\end{align*}
It is then clear that we can choose a sufficiently small $\alpha > 0$ such that
\begin{equation*}
\P(e(S, S^c) \leq \alpha s)  \leq c_3 s^{5/2}\mathrm{e}^{-\frac{s}{10} H(\frac{s}D)D} \leq c_3 s^{5/2}\mathrm{e}^{-\frac{s}{10} \log(D/s)}~.
\end{equation*}
Combined with the fact that $D = (4+o(1))\epsilon^3 n\to \infty$ by \eqref{eq-kernel-size}, and summing
over the possible values of $s$, we deduce that
$$\P(\exists S \subset \K,~ d_\K(S) \leq |E(\K)|,~ e(S, S^c) < \alpha d_\K(S)) = o(1)~,$$
as required.
\end{proof}

\begin{proof}[\emph{\textbf{Proof of Proposition \ref{prop-2-core-structure}}}]
The required statement on the typical number of vertices and edges in the kernel, $|\K|$ and $|E(\K)|$ resp., has
already been established in Lemma~\ref{lem-kernel-expander} (along with the expansion properties of the kernel).
It remains to show that, w.h.p., $N_2$ and $N'_2$ are both $(2+o(1))\epsilon^2 n$, where $N_2$ is the number of degree-2
vertices in the 2-core, and $N'_2$ equals $N_2$ minus the number of vertices that belong to disjoint cycles in the 2-core.

The main issue left is to distinguish between $H$ and $H'$ (the graph before and after removing its disjoint cycles).

Let $p_2(x)$ be the probability that a $\Po(x)$ variable equals $2$:
$$p_2(x) = \mathrm{e}^{-x}\frac{x^2}{2}~.$$
By Definition~\ref{def-poisson-conf-model} and our assumption that $\Lambda =  (2 + o(1)) \epsilon$,
$$N_2\sim \Bin(n, p_2(\Lambda))~,\mbox{ and }~ p_2(\Lambda_C) = (2+ o(1)) \epsilon^2\,.$$
 Hence, by a standard concentration argument (using the fact that $\epsilon^2 n\to\infty$) we deduce that $N_2 = (2+o(1)) \epsilon^2 n$ w.h.p. Assume therefore that this is the case (i.e., there are $(2+o(1))\epsilon^2 n$ vertices of degree 2 in $H$), and that
$|E(\K)| = \big(2 + o(1)\big)\epsilon^3 n$.

We next consider the disjoint cycles in $H$. For a given degree-2 vertex $v$ in $H$, let $A_{v,k}$ denote the
event that $v$ belongs to such a disjoint cycle whose length is $k$, and further let $A_v=\cup_{k=1}^{N_2} A_{v,k}$.
To form a disjoint cycle of length $k$ via the configuration model, we must repeatedly match points of degree-2 vertices, that is
(noticing that $2|E(\K)|$ counts the total number of points to be matched via the configuration model):
\begin{align*}
  \P(A_{v,k}) &=
\frac{1}{2N_2 +2|E(\K)|-2k+1} \prod_{j=0}^{k-2} \frac{2 N_2 - 2j-1}{2 N_2 + 2 |E(\K)| - 2j-1} \\
&=\frac{1}{2N_2 +2|E(\K)|-1} \prod_{j=0}^{k-2} \frac{2 N_2 - 2j-1}{2 N_2 + 2 |E(\K)| - 2j-3} \\
&\leq \frac{1+o(1)}{4\epsilon^2 n} \left(1-(1 + o(1))\epsilon\right)^{k-1}~,
\end{align*}
since the terms in the above product over $j$ formed a decreasing sequence. Summing over the values of $k$,
\begin{align}
\P (A_v) &= \sum_{k=1}^{N_2} \P(A_{v,k}) \leq
\frac{1+o(1)}{4\epsilon^2 n} \sum_{k=1}^{N_2} \left(1-(1+ o(1))\epsilon\right)^{k-1} \nonumber\\
&\leq  \frac{1+o(1)}{4\epsilon^2 n} \cdot \frac{1+o(1)}{\epsilon} = \frac{1+ o(1)}{4\epsilon^3 n} = o(1)\label{eq-prob-separate-cycle}~.
\end{align}
It then follows that w.h.p.\ $N'_2 = (1 - o(1)) N_2 = (2+ o(1))\epsilon^2 n$, as required.
\end{proof}

\subsection{Contiguity of Poisson-cloning and Poisson-configuration}

A key part of showing the contiguity result is a counterpart for Theorem~\ref{thm-Lambda-C-upper}, which together implies that $\Lambda_C$ has a tight concentration window of order $1/\sqrt{\epsilon n}$.

\begin{theorem}\label{thm-Lambda-C-lower}[Lower bound on the window of $\Lambda_C$]
There exist some constant $c > 0$ such that for any $t =t(n) > 0$ and fixed $\delta > 0$,
\begin{equation}\label{eq-Lambda-C-interval}
\P\left(t \leq \Lambda_C \leq t + \frac{\delta}{\sqrt{\epsilon n}}\right) \leq c \delta~.\end{equation}
\end{theorem}
\begin{proof}
The results of the previous subsection imply that, for some suitably large constant $M = M(\delta)> 0$, $$\P\Big(|\Lambda_C - \lambda \theta_\lambda| \geq \frac{M}{\sqrt{\epsilon
n}}\Big) \leq \delta~. $$
It follows that \eqref{eq-Lambda-C-interval} holds trivially for any $c \geq 1$ when $|t - \lambda \theta_\lambda| \geq \frac{M+1}{\sqrt{\epsilon n}}$.
Therefore, we assume in what follows $|t - \lambda \theta_\lambda| <\frac{M+1}{\sqrt{\epsilon n}}$.
Denote by $A$ the event $\{t \leq \Lambda_C \leq t+ \delta/\sqrt{\epsilon n}\}$. Recall the fact that $|\TC[G]| \sim \Bin(n, p^+_{2}(\Lambda_C))$, where $\TC[G]$ is the 2-core of the Poisson-cloning model and $p^+_{2}(x)$ stands for the probability for a $\Po(x)$ variable to be at least 2. Standard analytical arguments give that $$ p^+_{2}\Big(t + \frac{\delta}{\sqrt{\epsilon M}}\Big) - p^+_{2}(t) \leq \frac{3\delta}{\sqrt{n/\epsilon}}~.$$
Now, an application of CLT implies that for some interval $B$ of length $4 \delta \sqrt{\epsilon n}$,
\begin{equation}\label{eq-P-A-B}\P\left(|\TC[G]| \in B\given A\right) = 1-o(1)~.\end{equation} Consider $\TC$, the 2-core of the giant component in the Poisson-cloning. Recalling \eqref{eq-prob-separate-cycle} (which, as the discussion before Proposition~\ref{prop-2-core-structure}, also applies to the Poisson-cloning model), we know that w.h.p.\ only an $\frac{w}{\epsilon^3 n}$ fraction of vertices in $\TC[G]$ will appear in disjoint cycles, where $w = (\epsilon^3 n)^{1/4}~(\to\infty)$.
Therefore, we have
$$|\TC| = |\TC[G]|\Big( 1 -
\tfrac{w}{\epsilon^3 n}\Big) = |\TC[G]| -
O\big(\tfrac{w}{\epsilon}\big) = |\TC[G]| - o(\sqrt{\epsilon n})~.$$ Together with \eqref{eq-P-A-B}, we conclude that there exists an interval $B'$ with length $5 \delta \sqrt{\epsilon n}$ such that
\begin{equation}\label{eq-P-A-B'}\P\left(|\TC| \in B'\given A\right) = 1-o(1)~.\end{equation} Now, for the 2-core $H$ of the giant component of the Erd\H{o}s-R\'enyi graph $\cG(n,p)$, it is known (see \cite{PW}*{Theorem 6}, reformulated here in Theorem~\ref{thm-PW-local-limit})
that $|H|$ is in the limit Gaussian with variance $(12+o(1))\epsilon n$. Therefore,
\[\P(|H| \in B') \leq 5 \delta~.\]
Combining this with contiguity of Poisson-cloning and $\cG(n,p)$ (as stated in Theorem~\ref{thm-poisson-ER}), we obtain that for some constant $c_0 >0$ $$\P(|\TC| \in B') \leq 5 c_0 \delta~.$$ The proof is completed by choosing $c = 5 c_0 + 1$ and applying \eqref{eq-P-A-B'}.
\end{proof}

Using the above estimate for $\Lambda_C$, we are now able to conclude the main result of this section, which reduces the $2$-core of Poisson-cloning to the graph generated by the Poisson-configuration model.

\begin{proof}[\textbf{\emph{Proof of Theorem~\ref{thm-Lambda-contiguity}}}]
Recall that $H$ is the random graph generated by the Poisson-configuration model, and $\TC[G]$ is the 2-core of Poisson-cloning.

Let $\delta > 0$, and set
\[B = (\lambda \theta_\lambda - M\sqrt{\epsilon n},~ \lambda \theta_\lambda + M\sqrt{\epsilon n})\,,\] where $M=M(\delta)$ is a sufficiently large constant such that $\P(\Lambda_C \in B) \geq 1 - \delta$. Further define
\[ D \deq \{x: \P(H\in \mathcal{A} \mid \Lambda = x) \geq \delta\}\,.\]
 Since $\P(H\in \mathcal{A}) = o(1)$, we obtain that $\P(\Lambda \in D) \to 0$, and consequently \[\P(\Lambda \in B \cap D)=o(1)\,.\] Recalling that $\Lambda \sim \mathcal{N}\left(\lambda \theta_\lambda, 1/(\epsilon n)\right)$, we deduce that $\mathcal{L}(B \cap D)\sqrt{\epsilon n} \to 0$, where $\mathcal{L}(\cdot)$ stands for the Lebesgue measure on $\R$. At this point, Theorem~\ref{thm-Lambda-C-lower} gives that $\P(\Lambda_C \in B \cap D) \to 0$. Recalling
 that $\TC[G]$ and $H$ are generated by the same scheme (and hence have the same distribution) given the event $\Lambda_C = \Lambda$,
 we obtain that \[\P(\TC[G] \in \mathcal{A}) \leq 2\delta + o(1)~,\]
 as required.
\end{proof}

\section{Constructing the 2-core of the random graph}\label{sec:poisson-geometric}

In the previous section, we have shown that the 2-core of Poisson-cloning is contiguous to a simpler model, which we called
the Poisson-configuration model (see Definition~\ref{def-poisson-conf-model}).
The goal of this section is to reduce the Poisson-configuration model to the following, where here
 and in what follows, $\mu$ is defined to be the conjugate of $\lambda = 1+ \epsilon$. That is to say, $\mu < 1$ and
\begin{equation}
  \label{eq-mu-def}
  \mu\mathrm{e}^{-\mu} = \lambda \mathrm{e}^{-\lambda}~.
\end{equation}

\begin{definition}[Poisson-geometric model for $n$ and $p=\frac{1+\epsilon}n$]\label{def-poisson-geo-model}\mbox{}
\begin{enumerate}[(1)]
\item Let $\Lambda\sim \mathcal{N}\left(1 + \epsilon - \mu, \frac1{\epsilon n}\right)$ and assign an independent $\Po(\Lambda)$ variable $D_u$ to each vertex $u$. Let $N_k = \#\{u : D_u = k\}$ and $N = \sum_{k\geq 3}N_k$.
\item Construct a random graph $\K$ on $N$ vertices, uniformly chosen over all graphs with $N_k$ degree-$k$ vertices for $k\geq 3$ (if $\sum_{k\geq 3}k N_k$ is odd, choose a vertex $u$ with $D_u=k \geq 3$ with probability proportional to $k$, and give it $k-1$ half-edges and a self-loop).
\item Replace the edges of $\K$ by paths of length i.i.d.\ $\Geom(1-\mu)$.
\end{enumerate}
\end{definition}

\begin{theorem}\label{thm-Poisson-contiguity}
Let $H$ be generated by the Poisson-configuration model w.r.t.\ $n$ and
$p=\frac{1+\epsilon}n$, where $\epsilon\to 0$ and $\epsilon^3 n \to \infty$.
Let $\tilde{H}$ be generated by the Poisson-geometric model corresponding to $n,p$. Then for any set of graphs $\mathcal{A}$ such that $\P(\tilde{H} \in
\mathcal{A}) \to 0$, we have $\P(H \in \mathcal{A}) \to 0$.
\end{theorem}

Clearly, both models have the same kernel, and they only differ in the way this kernel is thereafter expanded to form the entire graph (replacing edges by paths). To prove the above statement, we need to estimate the distribution of the total number of edges in each of the models; we will show that they are in fact contiguous.

\subsection{Edge distribution in the Poisson-configuration model}

\begin{lemma}\label{eq-P(K|Lambda)-ratio}
Let $N_k$ denote the number of degree-$k$ vertices in the Poisson-configuration model,
and set $\Lambda_0 = \lambda - \mu$.
For any fixed $M > 0$ there exist some $c_1,c_2>0$ such that the following holds:
If $n_3,n_4,\ldots$ satisfy
\begin{align*}
\begin{array}{rl}
\Big|n\left(1-\mathrm{e}^{-\Lambda_0}(1+\Lambda_0 + \tfrac{\Lambda_0^2}2)\right) - \sum_{k\geq3}n_k \Big| &\leq M\sqrt{\epsilon^3 n}\,,\\
\Big|n\Lambda_0\left(1-\mathrm{e}^{-\Lambda_0}(1 + \Lambda_0)\right) - \sum_{k\geq3}k n_k \Big| &\leq M\sqrt{\epsilon^3 n}\,
\end{array}
\end{align*}
and $x$ satisfies $|x - \Lambda_0| \leq \frac{M}{\sqrt{\epsilon n}}$ then
\[
c_1 \leq \frac{\P\left(N_k = n_k~\mbox{ for all }k\geq 3 \given \Lambda = x\right)}
{\P\left(N_k = n_k~\mbox{ for all }k\geq 3 \given \Lambda = \Lambda_0\right)} \leq c_2\,.
\]
\end{lemma}
\begin{proof}
Throughout the proof of the lemma, the implicit constants in the $O(\cdot)$ notation depend on $M$.

Write $m = \sum_{k \geq 3}n_k$ and $r=\sum_{k\geq 3}k n_k$, and let
$A=A(n_3,n_4,\ldots)$ denote the event $\{N_k = n_k\mbox{ for all }k\geq 3\}$.
 As usual, we use the abbreviations
 \[p_k(x) = \P(\Po(x)=k) = \mathrm{e}^{-x}x^k/k!~,~\mbox{ and }p^-_k = \P(\Po(x) \leq k)\,.\] It follows that
\begin{align*}
\frac{\P\left(A \given \Lambda = x\right)}
{\P\left(A \given \Lambda = \Lambda_0\right)} &=
\left(\frac{p^-_2(x)}{p^-_2(\Lambda_0)}\right)^{n-m} \prod_k \left(\frac{p_k(x)}{p_k(\Lambda_0)}\right)^{n_k} \\
& = \mathrm{e}^{-n(x-\Lambda_0)}\left(\frac{1+x+\frac{x^2}2}{1+\Lambda_0+\frac{\Lambda_0^2}2}\right)^{n-m} \left(\frac{x}{\Lambda_0}\right)^r~,
\end{align*}
and so
\begin{align*}
\log\frac{\P\left(A \given \Lambda = x\right)}
{\P\left(A \given \Lambda = \Lambda_0\right)} &=
n(\Lambda_0-x) + (n-m)\log\bigg(\frac{1+x+\tfrac{x^2}2}{1+\Lambda_0+\frac{\Lambda_0^2}2}\bigg)+ r\log\frac{x}{\Lambda_0}\,.
\end{align*}
Using Taylor's expansion and recalling that $x -\Lambda_0 = O(1/\sqrt{\epsilon n}) = o(\Lambda_0)$,
\begin{align*}
\log\bigg(\frac{1+x+\tfrac{x^2}2}{1+\Lambda_0+\frac{\Lambda_0^2}2}\bigg)
 &= \frac{1+\Lambda_0}{1+\Lambda_0+\frac{\Lambda_0^2}2}(x-\Lambda_0)  -(\Lambda_0-o(1))(x-\Lambda_0)^2 \\
 &=\frac{1+\Lambda_0}{1+\Lambda_0+\frac{\Lambda_0^2}2}(x-\Lambda_0)  +O(1/n)\,,
\end{align*}
and we
deduce that
\begin{align*}
\log\frac{\P\left(A \given \Lambda = x\right)}
{\P\left(A \given \Lambda = \Lambda_0\right)} &=
n(\Lambda_0-x) +  (n-m) \frac{1+\Lambda_0}{1+\Lambda_0+\frac{\Lambda_0^2}2}(x-\Lambda_0)
- O(1) \\
&+ r \frac{x - \Lambda_0}{\Lambda_0} - O\left(r/ \epsilon^3 n \right)\,.
\end{align*}
Our assumptions on $m,r$ now yield that
\begin{align*}
\log\frac{\P\left(A \given \Lambda = x\right)}
{\P\left(A \given \Lambda = \Lambda_0\right)} &=
n(\Lambda_0-x) +  n \mathrm{e}^{-\Lambda_0} (1+\Lambda_0)(x-\Lambda_0)\\
&+ n \left(1-\mathrm{e}^{-\Lambda_0}(1+\Lambda_0)\right)(x-\Lambda_0) + O(1) = O(1)\,,
\end{align*}
completing the proof.
\end{proof}
Fix $M > 0$, and let $\cB_M$ denote the following set of ``good'' kernels:
\begin{align}\label{eq-good-kernels}
\cB_M \deq \left\{\K :\begin{array}{c}
\Big||\K| - n\left(1-\mathrm{e}^{-\Lambda_0}(1+\Lambda_0 + \tfrac{\Lambda_0^2}2)\right) \Big| \leq M\sqrt{\epsilon^3 n}\\
\Big||E(\K)| - \tfrac12 n\Lambda_0\left(1-\mathrm{e}^{-\Lambda_0}(1 + \Lambda_0)\right) \Big| \leq M\sqrt{\epsilon^3 n}
\end{array}\right\}\,.
\end{align}
Let $f_{\Lambda}(\cdot \mid \cdot)$ denote the density function of $\Lambda$ given the kernel $\K$
(or equivalently, given its degree sequence). By applying Bayes' formula, the above lemma gives that
\[ \frac{f_{\Lambda}(x \mid \K)}{f_{\Lambda}(\Lambda_0 \mid \K)} = \Theta(1)
\]
for all $\K \in \cB_M$ and $x$ in the interval $I_M = [\Lambda_0 - \frac{M}{\sqrt{\epsilon n}},\Lambda_0 + \frac{M}{\sqrt{\epsilon n}}]$. Clearly, by volume considerations, this implies that
for some $c = c(M) >0$ we have
\begin{equation}\label{eq-Lambda-given-K}
f_{\Lambda}(x \mid \K) \leq c\sqrt{\epsilon n}\quad\mbox{ for all $x \in I_M$ and $\K \in \cB_M$.}
\end{equation}

\begin{lemma}\label{lem-poisson-conf-upper}
Define $M>0$, $I_M$ and $\mathcal{B}_M$ as above. Let $H$ be generated by the
Poisson-configuration model. There exists some constant $c = c(M) > 0$ so that
for any $\K \in \mathcal{B}_M$ and $s$ with $\big|s - \frac{n}2\left(\Lambda_0 - \mathrm{e}^{-\Lambda_0}\Lambda_0 \right)\big| \leq M\sqrt{\epsilon n}$,
$$\P\left(|E(H)| = s \,,\, \Lambda \in I_M \given \ker(H) = \K\right) \leq \frac{c}{\sqrt{\epsilon n}}\,.$$
\end{lemma}
\begin{proof}
Let $x \in I_M$ and $\K \in \cB_M$, and write $m = |\K|$ and $r=|E(\K)|$ for the number of vertices and the edges in the kernel respectively. We will first estimate $\P(|E(H)| = s \given \Lambda = x\,,\, \ker(H) = \K)$, and the required inequality will then
readily follow from an integration over $x \in I_M$.

Note that, given $\Lambda = x$ and $\ker(H)=\K$, the number of edges in $H$ is the $r$ edges of $\K$ plus an added edge for each degree $2$ variable out of the $n-m$ variables (i.i.d.\ $\Po(x)$) that have $\{ u : D_u \leq 2\}$. That is, in this case
\[
|E(H)| \sim r+\Bin\left( n-m, \frac{x^2/2}{1 + x + x^2/2}\right)\,,
\]
and therefore,
\begin{align*}
\P&\left(|E(H)| = s \given \Lambda = x\,,\, \ker(H) = \K\right)\\
 &=\binom{n-m}{s-r}\left(\frac{x^2/2}{1+x+x^2/2}\right)^{s-r}\left(\frac{1+x}{1+x+x^2/2}\right)^{n-m-(s-r)}\,.
\end{align*}
Write
\[q_0 = \frac{\Lambda_0^2/2}{1+ \Lambda_0 +\Lambda_0^2}\] and define
$$q = q(x) = \frac{x^2/2}{1+ x + x^2/2}\qquad \mbox{ and } \qquad t = \frac{s-r}{n-m}~.$$
Since $x\in I_M$ and $\K \in \mathcal{B}_M$, we have
\begin{equation}
  \label{eq-q-t-values}
  q = q_0 + O(\sqrt{\epsilon/ n})~, \mbox{ and } t = q_0 + O(\sqrt{\epsilon /n})~.
\end{equation}
Using Stirling's formula, we obtain that
\begin{align*}
\P&\left(|E(H)| = s \given \Lambda = x\,,\, \ker(H) = \K\right) \\
&\leq \bigg(\frac{(1+o(1))(n-m)}{2\pi(s-r)\left(n-m-(s-r)\right)}\bigg)^{1/2} \left(\frac{q}{t}\right)^{s-r} \left(\frac{1-q}{1-t}\right)^{n-m-s+r}\\
&\leq \frac{1}{\sqrt{s - r}} \mathrm{e}^{(n-m)g_t(x)}~,
\end{align*}
where $g_t(q)$ is given by
\begin{align*}
g_t(q) \deq - t \log t  - (1-t) \log (1-t) + t \log q + (1-t)\log(1-q)~.
\end{align*}
It is easy to verify that
$$ g_t'(q) = -\frac{1-t}{1-q} + \frac{t}{q}\,,$$
and so for any $q,t$ satisfying \eqref{eq-q-t-values} we have
\begin{align*}
g_t(t) &= 0\,,\qquad g'_t(t) = 0\,,\\
g''_t (q)&= -\frac{t}{q^2} - \frac{1 - t}{(1-q)^2} = -\frac{1+o(1)}{q_0}\,.
\end{align*}
Thus, for any large $n$ (absorbing the $o(1$)-term in the constant) we have
\begin{equation}
  \label{eq-gt(q)-leq}
  g_t(q) \leq -\frac{1}{3 q_0} (q-t)^2~.
\end{equation}

Clearly, the function $q(x) = \frac{x^2/2}{1 + x + x^2/2}$ satisfies
$q'(x) = \frac{x(2+x)}{2(1+x+x^2/2)^2}$, and in particular $q$ is strictly monotone increasing from $0$ to $1$
for $x\in[0,\infty)$. Thus, there exists a unique $x_t>0$ such that $q(x_t) = t$.
Noticing that for all $x = (1+o(1))\Lambda_0$ we have $q'(x) = (1+o(1))\Lambda_0 = (2+o(1))\epsilon$,
it follows that for any $M_1 > 0$ one can choose $M_2 > 0$ such that
\[
\left[\Lambda_0 - M_1\sqrt{\epsilon/ n},\Lambda_0 + M_1\sqrt{\epsilon /n}\right] \subset q\left(\left[\Lambda_0 - \tfrac{M_2}{\sqrt{\epsilon n}},\Lambda_0 + \tfrac{M_2}{\sqrt{\epsilon n}}\right]\right)\,,
\]
and in particular, $x_t = \Lambda_0 + O(1/\sqrt{\epsilon n})$.
We can now apply the Mean Value Theorem to $q$ in \eqref{eq-gt(q)-leq} and obtain that
$$g_t(q) \leq -\frac{1}{3q_0} (q(x) - q(x_t))^2 = -\frac{((1+o(1))\Lambda_0)^2}{3q_0}(x - x_t)^2 \leq -\frac{3}{5}(x - x_t)^2\,,$$
where the last inequality holds for any sufficiently large $n$.
Altogether, absorbing the change from $(n-m)$ to $n$ in the constant, we conclude that
\begin{align}\label{eq-prob-conditional}
\P&(|E(H)| = s \mid \Lambda = x,~ \ker(H) = \K)\nonumber \\
&\leq \frac{1}{\sqrt{s - r}} \mathrm{e}^{-\frac{1}{2}n(x - x_t)^2} \leq \frac{1}{\sqrt{\epsilon^2
n}}\mathrm{e}^{-\frac{1}{2}n(x - x_t)^2}~.
\end{align}

It remains to integrate the above conditional probability over $x \in I_M$. Combining \eqref{eq-Lambda-given-K} and
\eqref{eq-prob-conditional}, we obtain that for some constant $c > 0$
\begin{align*}
\P&\left(|E(H)| = s~,~ \Lambda\in I_M \given \ker(H) = \K\right)\\
& \leq \int_{I_M} \P\left(|E(H)| = s \given \Lambda = x~,~ \ker(H) = \K\right) f_\Lambda(x | \K)\, dx\\
&\leq \int_{I_M} \frac{1}{\sqrt{\epsilon^2 n}} \mathrm{e}^{-\frac{1}{2}n(x - x_t)^2} c\sqrt{\epsilon n} \leq \frac{c}{\sqrt{\epsilon n}} \int_{-\infty}^{\infty} \mathrm{e}^{-\frac{1}{2} y^2} dy =
\frac{c \sqrt{2\pi}}{\sqrt{\epsilon n}},
\end{align*}
as required.
\end{proof}

\subsection{Edge distribution in the Poisson-geometric model}

\begin{lemma}\label{lem-poisson-geo-lower}
Let $M>0$ and $\mathcal{B}_M$ be as in \eqref{eq-good-kernels}. Let $\tilde{H}$ be generated by the Poisson-geometric model. There exists some constant $c = c(M) > 0$ so that
for any $\K \in \mathcal{B}_M$ and $s$ with $\big|s - \frac{n}2(\lambda-\mu)\left(1-\frac{\mu}{\lambda}\right)\big| \leq M\sqrt{\epsilon n}$,
$$\P\left(|E(\tilde{H})| = s \given \ker(\tilde{H}) = \K\right) \geq \frac{c}{\sqrt{\epsilon n}}\,.$$
\end{lemma}
\begin{proof}
By definition, given that $\ker(\tilde{H})=\K$, the variable $|E(\tilde{H})$ is the sum of $|E(\K)|$ i.i.d.\ geometric variables
with mean $1/(1-\mu)$.

Denote by $r$ the number of edges in the kernel $\K$, and let $s$ be a candidate for the number of edges in the expanded 2-core $\tilde{H}$. As stated in the lemma (recall definition \eqref{eq-good-kernels}), we are interested in the following range for $r,s$:
\begin{align*}
r &= \frac{n}{2}(\lambda - \mu) (1 - \tfrac{\mu}{\lambda})(1 - \mu) +
c_1 \sqrt{\epsilon^3 n}, &\quad (|c_1| \leq M)~,\\
s & = \frac{n}{2}(\lambda - \mu) (1 - \tfrac{\mu}{\lambda}) + c_2 \sqrt{\epsilon n}, &\quad (|c_2| \leq M)~.
\end{align*}
In this case, we have that
\begin{align}\label{eq-ratio-n1-n2}
\frac{s - 1}{r -1} 
&= \frac{1}{1-\mu} + \frac{c_2\sqrt{\epsilon n} - c_1 \sqrt{\epsilon^3
n}/(1-\mu) - 1}{\frac{n}{2}(\lambda - \mu) (1 - \tfrac{\mu}{\lambda})(1 - \mu) + c_1 \sqrt{\epsilon^3 n}-1}
= \frac{1+\xi}{1-\mu}~, \end{align} where $\xi \deq\frac{1+o(1)}{2}\frac{c_2 - c_1}{\sqrt{\epsilon^3 n}}$.
Let $X_i$ be independent geometric random variables with mean $\frac{1}{1-\mu}$, i.e., $\P(X_i=k) = \mu^{k-1}(1-\mu)$ for $k=1,2,\ldots$; further set $S_k = \sum_{i=1}^k (X_i - 1)$. Since $S_k$ follows a negative binomial distribution, \begin{align*} \P(S_{r} = s - r) = \binom{s -1}{r -1} (1-\mu)^{r} \mu^{s - r}~.
\end{align*}
Using Stirling's formula, we get that for some constant $c_3>0$ \begin{align*} \P(S_{r} = s - r) \geq \frac{c_3}{\sqrt{r}}\Big(\frac{s - 1}{r -1}\Big)^{r - 1} \Big(\frac{s - 1}{s - r}\Big)^{s - r}(1-\mu)^{r} \mu^{s - r}~.
\end{align*}
Substituting \eqref{eq-ratio-n1-n2} in the above, and using the fact that
\[\frac{s-1}{s-r} = \frac{(1+\xi)/(1-\mu)}{(1+\xi)/(1-\mu) - 1} = \frac{1+\xi}{\xi+\mu}~,\]
we obtain that \begin{align*} \P(S_{r} = s - r) &\geq \frac{c_3}{\sqrt{r}}\Big(\frac{1
+\xi}{1-\mu}\Big)^{r - 1} \Big(\frac{1 + \xi}{\xi+\mu}\Big)^{s -
r}(1-\mu)^{r} \mu^{s - r}\\
&= \frac{c_3 (1 - \mu)}{\sqrt{r}}(1 + \xi)^{s -1} \Big( \frac{\mu}{\xi
+ \mu}\Big)^{s - r}\\
&= \frac{c_3 (1 - \mu)}{\sqrt{r}}\exp\Big(g(\xi) \frac{r -
1}{1 - \mu}\Big)~,
\end{align*}
where \[ g(x) \deq (1 + x)\log (1+x) + (x + \mu)\log\big(\frac{\mu}{x
+ \mu}\big)~.\] Clearly, we have that $g(0) = 0$ and a standard
calculation yields that
$$g'(x) = \log (1 + x) + \log \Big(\frac{\mu}{x+\mu}\Big) \mbox{ and } g''(x) = \frac{\mu -1}{(x+\mu)(1+x)}~.$$ In particular, we have that $g'(0) = 0$ and $|g''(x)| \leq 2 (1 - \mu)$ when $|x| \leq |\xi|$, where $\xi$ is defined as above (recall that $\xi=o(1)$). Therefore, for all such $x$ we have $|g(x)| \leq 2(1-\mu)x^2$, and altogether,
\begin{align*} \P(S_{r} = s - r) &\geq \frac{c_3 (1 - \mu)}{\sqrt{r}}
\mathrm{e}^{-2(r - 1) |\xi|^2}~.
\end{align*}
Since
\[ r |\xi|^2 = \frac{1+o(1)}2(c_2-c_1)^2 \leq (2+o(1))M^2~,\]
we deduce that for some $c'_3 > 0$,
\begin{align*}\P(S_{r} = s - r) &\geq c'_3 / \sqrt{\epsilon n}\,,
\end{align*}
completing the proof.
\end{proof}

\subsection{Contiguity of the two models}
We are now ready to prove the main result of this section, Theorem~\ref{thm-Poisson-contiguity}, which reduces the Poisson-configuration model to the Poisson-geometric model.
\begin{proof}
For some constant $M > 0$ to be specified later, define the event
\begin{align*}
A_M \deq \left\{~\Lambda \in I_M,~ \ker(H) \in \cB_M,~ \left| |E(H)| - \tfrac{n}{2}(\lambda - \mu) (1 - \tfrac{\mu}{\lambda})\right| \leq M \sqrt{\epsilon n}~
   \right\}\,.
\end{align*}
Fix $\delta > 0$. We claim that for a sufficiently large
$M = M(\delta)$ we have $\P(A_M) \geq 1- \delta$.
To see this, note the following:
\begin{enumerate}[\indent 1.]
  \item In the Poisson-configuration model, $\Lambda \sim \mathcal{N}\left(\Lambda_0 , \frac{1}{\epsilon n}\right)$, and $I_M$ includes at least $M$ standard deviations about its mean.
      \item Each of the variables
      $|\K |$ and $E(\K)$ is a sum of i.i.d.\ random variables with variance $O(\epsilon^3 n)$ and mean as specified in the definition of $\cB_M$, hence their concentration follows from CLT.
      \item Finally, $E(H)$ is again a sum of i.i.d.\ variables and has variance $O(\epsilon n)$, only here we must subtract the vertices that comprise disjoint cycles. By \eqref{eq-prob-separate-cycle} and the estimate in Proposition~\ref{prop-2-core-structure} on the size of the 2-core in the Poisson-configuration model, the number of such vertices is $O(1/\epsilon)$ w.h.p. Compared to the standard deviation of $O(\sqrt{\epsilon n})$, this amounts to a negligible error, as $\epsilon^3 n\to \infty$.
\end{enumerate}
Given an integer $s$ and a kernel $\K$, let $\mathcal{D}_{s,\K}$ denote every possible $2$-core with $s$ edges and kernel $\K$.
Crucially, the distribution of the Poisson-configuration model given $E(H)=s$ and $\ker(H)=\K$ is uniform over $\mathcal{D}_{s,\K}$, and so is the Poisson-geometric model given $E(\tilde{H})=s$ and $\ker(\tilde{H})=\K$.
Therefore, for any graph $D \in \mathcal{D}_{s,\K}$,
\begin{align*}
\frac{\P(H = D\given \ker(H)=\K )}{\P(\tilde{H} = D\given \ker(\tilde{H})=\K )}
= \frac{\P (|E(H)| = s \given \ker(H)=\K )}{\P (|E(\tilde{H})| = s \given \ker(\tilde{H})=\K)}~.
\end{align*}
Combining Lemmas~\ref{lem-poisson-conf-upper},\ref{lem-poisson-geo-lower} we get that for some $c=c(M) > 0$,
\[
\frac{\P (|E(H)| = s~,~A_M \given \ker(H)=\K )}{\P (|E(\tilde{H})| = s \given \ker(\tilde{H})=\K)}
\leq c~.\]
Recalling that $\P(A_M) \geq 1 -\delta$ and letting $\delta \to 0$, we deduce that for any family of graphs $\mathcal{A}$,
if $\P(\tilde{H} \in \mathcal{A}) \to 0$ then also $\P(H \in \mathcal{A}) \to 0$.
\end{proof}

\section{Constructing the giant component}\label{sec:structure-gc}

Throughout the section, let $p = (1+\epsilon)/n$, where $\epsilon \to 0$ and $\epsilon^3 n \to \infty$ with $n$, and let
$G$ be a random graph $G \sim \cG(n,p)$. We begin by analyzing
the ``bushes'', i.e., the trees that are attached to $\TC[G]$, the 2-core of $G$.

As before, $\mu < 1$ is defined to be the conjugate of $\lambda = 1+ \epsilon$ (see \eqref{eq-mu-def}).
Since $\epsilon \to 0$, we can infer from a standard Taylor expansion that
\begin{equation}
  \label{eq-mu-taylor}
\mu = 1 - \epsilon + \tfrac23 \epsilon^2 +O(\epsilon^3)~.
\end{equation}

%

\begin{proof}[\emph{\textbf{Proof of Theorem~\ref{mainthm-struct-gen}}}]
In what follows, we use the abbreviation PGW($\mu$)-tree for a Poisson($\mu$)-Galton-Watson tree.
Let $\hGC$ denote the graph obtained as follows:
\begin{itemize}
  \item Let $H$ be a copy of $\TC$ (the 2-core of the giant component of $G$).
  \item For each $v \in H$, attach an independent PGW($\mu$) tree rooted at $v$.
\end{itemize}
By this definition, $\GC$ and $\hGC$ have exactly the same 2-core $H$. For simplicity, we will refer directly to $H$ as the 2-core of the model, whenever the context of either $\GC$ or $\hGC$ is clear. We first establish the contiguity of $\GC$ and $\hGC$.

Define the bushes of $\GC$ as follows:
$$T_u \deq \{v\in \GC: v \mbox{ is connected to } u \mbox{ in } \GC \setminus H\}~\mbox{ for $u\in H$}\,.$$
Clearly, each $T_u$ is a tree as it is connected and does not contain any cycles (its vertices were not included in the 2-core).
To conclude, we go from $H$ to $\GC$ by attaching a tree $T_u$ to each vertex $u\in H$ (while identifying the root of $T_u$ with $u$).
Analogously, let $\{\tilde{T}_u\}_{u\in H}$ be the corresponding bushes in $\hGC$.

We next introduce notations for the labeled and unlabeled trees as well as their distributions. For $t \in \N$, let $\mathcal{R}_t$ be the set of all labeled rooted trees on the vertex set $[t]$, and let $U_t$ be chosen uniformly at random from $\mathcal{R}_t$. For $T \in \mathcal{R}_t$ and a bijection $\phi$ on $[t]$, let $\phi(T)$ be the tree obtained by relabeling the vertices in $T$ according to $\phi$. Furthermore, define
$$T' \deq \{\phi(T): \phi \mbox{ is a bijection on }[t]\}$$
to be the corresponding rooted unlabeled tree.

Let $\{t_u : u \in H\}$ be some integers. Conditioning on the event
$$\{\,|T_u| = t_u \mbox{ for all } u\in H\,\}\,,$$
we know from the definition of $\cG(n,p)$ that $T_u$ is distributed independently and uniformly among all labeled trees of size $t_u$ rooted at $u$. In particular, in that case each $T'_u$ is independently distributed as $U'_{t_u}$ (the unlabeled counterparts of $T_u$ and $U_{t_u}$).

On the other hand, Aldous \cite{Aldous} (see also, e.g., \cite{AP}) observed that, if $T$ is a PGW-tree then $T'$ has the same distribution as $U'_t$ on the event $\{|T| = t\}$. Therefore, conditioning on the event
$$\{\,|\tilde{T}_u| = t_u \mbox{ for all } u\in H\,\}~,$$
we also get that $\tilde{T}'_k$ has the same distribution as $U'_{t_k}$.

We therefore turn to study the sizes of the bushes in $\GC$ and $\hGC$. Letting $\{t_u : u \in H\}$ be some integers and writing
$$N = \sum_{u \in H} t_u\,,$$
we claim that by definition of $\cG(n,p)$ every extension of the 2-core $H$ to the component $\GC$, using trees whose sizes sum up to $N$, has the same probability. To see this, fix $H$, and notice that the probability of obtaining a component with a 2-core is $H$ and an extension $X$ connecting it to $N-|H|$ additional vertices only depends on the number of edges in $H$ and $X$ (and the fact that this is a legal configuration, i.e., $H$ is a valid 2-core and $X$ is comprised of trees). Therefore, upon conditioning on $H$ the probabilities of the various extensions $X$ remain all equal.
Cayley's formula gives that there are $m^{m-1}$ labeled rooted trees on $m$ vertices, and so,
\begin{align}\label{eq-size-tree-distr}
\P&\left(|T_u| = t_u \mbox{ for all } u \in H \given H \right) \nonumber\\
&= \P \left(|\GC| = N \given H\right) \P \left(|T_u| = t_u \mbox{ for all } u \in H  \given H\;,\; |\GC| = N\right)\nonumber\\
&= \P \left(|\GC| = N \given H \right) \frac{1}{Z(N)}\frac{N!}{\prod_{u\in H} t_u !} \prod_{u\in H} t_u^{t_u -1}\nonumber\\
&= \P (|\GC| = N \given H)\frac{1}{Z'(N)} \prod_{u \in H} \Big[\frac{t_u^{t_u -1}}{\mu t_u!} (\mu \mathrm{e}^{-\mu})^{t_u}\Big],
\end{align}
where $Z(N)$ and $Z'(N)$ are the following normalizing constants
\begin{align*}
Z'(N) &= \sum_{\{r_u\}:\sum_{u\in H}r_u = N} \prod_{u \in H} \Big[\frac{r_u^{r_u -1}}{\mu r_u!} (\mu \mathrm{e}^{-\mu})^{r_u}\Big]~,\\
 Z(N) &= Z'(N) \mu^{N-|H|} \mathrm{e}^{-\mu N}~.
\end{align*}
Notice that the size of a Poisson$(\gamma)$-Galton-Watson tree $T$ follows a Borel($\gamma$) distribution (see, e.g., \cite{Pitman}), namely,
\begin{equation}\label{eq-PGW-size}
\P(|T| = t)=\frac{t^{t-1}}{\gamma t!}(\gamma \mathrm{e}^{-\gamma})^{t}~.
\end{equation}
Recalling that $\tilde{T}_u$ are independent PGW($\mu$)-trees, it follows that
$$Z'(N) = \sum_{\{r_u\}:\sum_{u\in H}r_u = N} \Big[\prod_{u\in H}\P(|\tilde{T}_u| = r_u) \Big]= \P\left(|\hGC| = N \given H\right)~.$$
Combining this with \eqref{eq-size-tree-distr} and \eqref{eq-PGW-size}, we obtain that
\begin{align}\label{eq-prob-ratio}
\frac{\P(|T_u| = t_u \mbox{ for all } u \in H \given H )}{\P(|\tilde{T}_u| = t_u \mbox{ for all } u \in H \given H )}
= \frac{\P (|\GC| = N \given H)}{\P (|\hGC| = N \given H)}~.
\end{align}
At this point, we wish to estimate the ratio in the right hand side above. To this end, we need the following result of \cite{PW}, which we restate in our setting of the near-critical regime.
\begin{theorem}[\cite{PW}*{Theorem 6}, reformulated]\label{thm-PW-local-limit}
Let $\mathbf{b}(\lambda) = \left(\begin{smallmatrix} b_1(\lambda) \\ b_2(\lambda) \\ b_3(\lambda)\end{smallmatrix}\right)$
where
\begin{align*}
b_1(\lambda) = (1 - \mu)\left(1 - \tfrac{\mu}\lambda\right)\,,\, b_2(\lambda) = \mu \left(1 - \tfrac{\mu}\lambda\right)\,,\,b_3(\lambda) = \tfrac12\left(1 - \tfrac\mu\lambda\right)(\lambda + \mu -2)\,.
\end{align*} There exist positive definite matrices $K_p,K_m$ satisfying
\begin{align*}
K_p &= \left(\begin{array}{ccc}
(12+o(1))\epsilon & 4+o(1)& (6+o(1))\epsilon^2\\
4+o(1) & (2+o(1))/\epsilon & (2+o(1))\epsilon\\
(6+o(1))\epsilon^2&(2+o(1))\epsilon& (\tfrac{10}3+o(1)) \epsilon^3
\end{array}
\right)~,\\
K_m &= K_p - 2\lambda \frac{d \mathbf{b}(\lambda)}{d\lambda} \cdot \frac{d \mathbf{b}(\lambda)^T}{d\lambda}~,
\quad \det(K_m) = (\tfrac83 +o(1))\epsilon^3~,
\end{align*}
and such that
\begin{enumerate}[(i)]
\item \label{item-PW-central-limit}$(|H|, |\GC| - |H|, |E(H)| - |H|)$ is in the limit Gaussian with a mean vector $n \mathbf{b}$ and a covariance matrix $n K_p$.
\item \label{item-PW-local-limit} If $A_m \deq K_m^{-1}$ and $B$ denotes the event that $|E(G)|=m$ for some $m= (1 + (1+o(1))\epsilon)\frac{n}{2}$, and there is a unique component of size between $\epsilon n$ and $4 \epsilon n$ and none larger, then
\begin{align}
\P& \left(|H| = n_1, |\GC| - |H| = n_2, |E(H)| - |H| = n_3 \given B\right)\nonumber \\
&= \frac{\sqrt{3} +o(1)}{8(\pi n \epsilon)^{3/2}} \exp\left(-\tfrac{1}{2}\mathbf{x}^T A_m \mathbf{x}\right)\,,
\end{align}
uniformly for all $(n_1, n_2, n_3) \in \N^3$ such that
$$(K_p(1, 1)^{-1/2} x_1,K_p(2, 2)^{-1/2} x_2, K_p(3, 3)^{-1/2} x_3 )$$
is bounded, where $\mathbf{x}^T=(x_1,x_2,x_3)$ is defined by
$$\mathbf{x}^T =\frac{1}{\sqrt{n}}(n_1 - b_1 n, n_2 - b_2 n, n_3 - b_3 n)~.$$
\end{enumerate}
\end{theorem}
Since $\epsilon^2 n\to\infty$, it is clear by CLT that w.h.p.\ the total number of edges in $G \sim \cG(n,p)$ is $(1+ (1+o(1))\epsilon)\frac{n}{2}$. Furthermore, by the results of \cite{Bollobas84} and \cite{Luczak90} (see also \cite{JLR}), w.h.p.\ our graph $G$ has a unique giant component of size
$$(1 + o(1))(1 - \mu/\lambda)n = (2+o(1))\epsilon n~.$$
Altogether, we deduce that the event $B$ happens w.h.p.; assume therefore that $B$ indeed occurs.
Define the event $Q$ by
\begin{align*}\mathcal{Q}_{M} &\deq \left\{(n_1, n_2, n_3)\in \N^3:  |x_1|\leq \sqrt{\epsilon} M~,~ |x_2|\leq \frac{M}{\sqrt{\epsilon}}~,~ |x_3|\leq \epsilon^{3/2} M \right\}~,\\
Q &\deq \big\{(|H|, |\GC|- |H|, |E(H)| - |H|) \in \mathcal{Q}_M\big\}~.
\end{align*}
By part~\eqref{item-PW-central-limit} of Theorem~\ref{thm-PW-local-limit}, for any fixed $\delta > 0$ there exists some $M > 0$ such that $\P(Q^c) < \delta$ for a sufficiently large $n$. Next, define
\begin{align*}
P_{\max} & = \max_{(n_1, n_2, n_3)\in \mathcal{Q}_M} \P \left(|H| = n_1,\; |\GC| - |H| = n_2,\; |E(H)| - |H| = n_3 \right)\,,\\
P_{\min} & = \min_{(n_1, n_2, n_3)\in \mathcal{Q}_M} \P \left(|H| = n_1,\; |\GC| - |H| = n_2,\; |E(H)| - |H| = n_3 \right)\,.
\end{align*}
It follows from part~\eqref{item-PW-local-limit} of Theorem~\ref{thm-PW-local-limit} that there exists some $c=c(M) > 0$ such that
\begin{equation}\label{eq-P-M-P-m}P_{\max} \leq c \cdot P_{\min}~,\end{equation}
when $n$ is sufficiently large.
Notice that by definition of $\mathbf{x}$,
\begin{equation*}\#\{n_2 \in \N: |x_2|\leq M/\sqrt{\epsilon} \} \geq M \sqrt{n/\epsilon}~.\end{equation*}
Combined with \eqref{eq-P-M-P-m}, it follows that for any $(n_1, n_2, n_3)\in \mathcal{Q}_M$ we have
\begin{equation}\label{eq-G-local}
\P\left( |\GC|= n_1 + n_2~,~ Q  \given |H| = n_1\right) \leq \frac{c}{M \sqrt{n/\epsilon}}~.
\end{equation}
With this estimate for $\P(|\GC|=N \given H)$, the numerator in the right-hand-side of \eqref{eq-prob-ratio}, it remains to estimate the denominator, $\P(|\hGC|=N \given H)$.

Recall that, given $H$, the quantity $|\hGC|$ is a sum of $|H|$ i.i.d.\ Borel($\mu$) random variables (each such
variable is the size of a PGW($\mu$)-tree). We would now like to derive a local central limit theorem
for $|\hGC|$. Unfortunately, each Borel($\mu$) variable $|T_u|$ has $\var|T_u|\asymp 1/\epsilon^3 \to\infty$, and standard versions of local CLT do not apply here. To bypass this obstacle, we use a different characterization of the tree-sizes $\{|T_u|:u \in H\}$.

It is well known that the total progeny in a branching process with offspring distribution $Z$ has the same law
as the hitting time from $1$ to $0$ of a one-dimensional random walk whose increments are i.i.d.\ variables distributed as
$Z-1$ (see, e.g., \cite{Spitzer}*{page 234}). Hence, the total size of $k$ i.i.d.\ such branching processes is exactly the hitting
time of this walk from $k$ to $0$. The following theorem of Otter \cite{Otter} characterizes this quantity:
\begin{theorem}[\cite{Otter}, see also \cite{vdHK}]\label{thm-hitting-time-rw}
Let $W_t$ be a random walk, whose steps $Y_i$ are i.i.d.\ random variables satisfying $Y_i \geq -1$. Then
$$ \P_k(\tau_0 = t) = \frac{k}{t} \P_k(W_t = 0)~.$$
\end{theorem}
In our setting, $Y_i \sim \Po(\mu)-1$, hence $\P_k(W_t = 0)=\P(S_t = t-k)$ where $S_t = \sum_{i=1}^t X_i$ and the $X_i$'s
are i.i.d.\ $\Po(\mu)$ variables. In light of the above theorem, it follows that for any integer $n_2$,
\begin{equation}
  \label{eq-tGc-size}
\P(|\hGC|= n_1+n_2 \given |H|=n_1) = \frac{n_1}{n_1+n_2}\P(S_{n_1+ n_2}=n_2)~.
\end{equation}
Since $S_m \sim \Po(m \mu)$ for any $m$, we have that
\begin{equation}
  \label{eq-S[n1+n2]}
  \P(S_{n_1+ n_2}=n_2) = \mathrm{e}^{-(n_1+n_2)\mu} \frac{\left((n_1+n_2)\mu\right)^{n_2}}{n_2!}~.
\end{equation}
Recalling the definition of $\mathcal{Q}_m$, we are interested in the following range for $n_1$ and $n_2$:
\begin{align*} n_1 &= (1-\mu)(1-\tfrac{\mu}{\lambda})n + c_1 \sqrt{\epsilon n}&(|c_1| \leq M)~,\\
n_2 &= \mu(1-\tfrac{\mu}{\lambda})n + c_2 \sqrt{n/\epsilon}&(|c_2| \leq M)~.
\end{align*}
In this case, we have
 \begin{align*}\frac{n_1+n_2}{n_2} &= \frac{(1-\frac{\mu}{\lambda})n + c_1\sqrt{\epsilon n} + c_2\sqrt{n/\epsilon}} {\mu(1-\frac{\mu}{\lambda})n+c_2\sqrt{n/\epsilon}} = \frac{1}{\mu} + \frac{c_1\sqrt{\epsilon n} +c_2(1-\frac1{\mu})\sqrt{n/\epsilon}}
{\mu(1-\frac{\mu}{\lambda})n + c_2\sqrt{n/\epsilon}} \\ &\deq \frac{1}{\mu} + \xi ~,\end{align*}
where $\xi = \xi(n) = (\tfrac12+o(1))(c_1-c_2)/\sqrt{\epsilon n}$.
Applying Stirling's formula to \eqref{eq-S[n1+n2]} and using the fact that
$1+x \geq \exp(x-x^2)$ for $x \geq 0$ gives
\begin{align*}
  \P(S_{n_1+ n_2}=n_2) &= \exp\left[\left(1-\frac{n_1+n_2}{n_2}\mu\right)n_2\right] \frac{1}{\sqrt{2\pi n_2}}\left(\frac{(n_1+n_2)\mu}{n_2}\right)^{n_2} \\
  &= \exp\left(-\xi\mu n_2\right) \frac{1}{\sqrt{2\pi n_2}}\left(1 + \xi \mu\right)^{n_2} \geq
  \frac{1}{\sqrt{2\pi n_2}} \mathrm{e}^{-\xi^2 \mu^2 n_2}~.
\end{align*}
Now, since $n_2 = (2+o(1))\epsilon n$ and
$$ \xi^2 \mu^2 n_2 = (1+o(1))\frac{(c_1-c_2)^2}{4\epsilon n} 2\epsilon n \leq (2+o(1))M^2~, $$
we conclude that for some constant $\delta' = \delta'(M)$,
\[
  \P(S_{n_1+ n_2}=n_2) \geq  \frac{\delta'}{\sqrt{n \epsilon}} ~.
\]
Recalling that $\frac{n_1}{n_1+n_2} = (1+o(1))\epsilon$, we can decrease $\delta'$
to absorb this $o(1)$ error-term for a sufficiently large $n$, and together with \eqref{eq-tGc-size} get
\begin{equation}\label{eq-tilde-G-local} \P(|\hGC|= n_1+n_2 \given |H|=n_1) \geq \delta'\frac{\epsilon}{\sqrt{n \epsilon}} = \frac{\delta'}{\sqrt{n/\epsilon}}~.
\end{equation}
Combining \eqref{eq-G-local} and \eqref{eq-tilde-G-local}, we obtain that when $n$ is sufficiently large,
 $$\frac{\P \left(|\GC| = N~,~Q \given |H|\right)}{\P \left(|\hGC| =  N \given |H|\right)} \leq \frac{c}{M \delta'}~.$$
By \eqref{eq-prob-ratio} (and recalling the fact that conditioned on $|T_i|$, the tree $T_i$ is uniformly distributed among all unlabeled trees of this size, and a similar statement holds for $\tilde{T}_i$), we conclude that for some $c' =c'(M)> 0$ and any unlabeled graph $A$
\begin{align}\label{eq-GC-hat-C}
\P (\GC =A\,,Q\,,B \mid H) \leq c'\,\P (\hGC = A \mid H)~.
\end{align}

We are now ready to conclude the proof of the main theorem. Let $\tGC$ be defined as in Theorem~\ref{mainthm-struct-gen}. For any set of simple graphs $\mathcal{A}$, define
\begin{equation}
  \label{eq-cH-def}
  \mathcal{H} = \left\{H : \P(\GC \in \mathcal{A}\, , Q\,, B\mid \TC = H) \geq (\P(\tGC \in \mathcal{A}))^{1/2}\right\}\,.
\end{equation}

Recall that by definition, $\tGC$ is produced by first constructing its 2-core (first two steps of the description), then attaching to each of its vertices independent PGW($\mu$)-trees. Hence, for any $H$, the graphs $\hGC$ and $\tGC$ have the same conditional distribution given $\TC[\hGC] = \TC[\tGC] = H$. It then follows from \eqref{eq-GC-hat-C},\eqref{eq-cH-def} that for some constant $c''>0$ and any $H \in \mathcal{H}$,
\[\P(\tGC \in \mathcal{A} \mid \TC[\tGC] = H) \geq c'' (\P(\tGC \in \mathcal{A}))^{1/2}\,.\]
Since
\[ \P(\tGC \in \mathcal{A}) \geq c'' (\P(\tGC \in \mathcal{A}))^{1/2} \P(\TC[\tGC] \in \mathcal{H}) \,,\]
the assumption that $\P(\tGC \in \mathcal{A}) \to 0$ now gives that $\P(\TC[\tGC] \in \mathcal{H}) \to 0$.

At this point, we combine all the contiguity results thus far to claim that, for any family of simple graphs $\mathcal{F}$,
\[ \P(\TC[\tGC] \in \mathcal{F}) = o(1) \mbox{ implies that }\P(\TC \in \mathcal{F})=o(1)\,.\]
Indeed, by definition, the 2-core of $\tGC$ is precisely the Poisson-geometric model, conditioned on the sum of the degrees ($\sum_u D_u \one_{D_u \geq 3}$) being even. Therefore, as $\mathcal{F}$ consists only of simple graphs,
clearly we may consider this model condition on the graph produced being simple, and in particular, that $\sum_u D_u \one_{D_u \geq 3}$ is even.
Applying Theorem~\ref{thm-Poisson-contiguity} (contiguity with Poisson-configuration), Theorem~\ref{thm-Lambda-contiguity} (contiguity with Poisson-cloning) and Theorem~\ref{thm-poisson-ER} (contiguity with Erd\H{o}s-R\'enyi graphs), in that order, now gives the above statement.

This fact and the arguments above now give that $\P(\TC \in \mathcal{H}) \to 0$. By the definition of $\mathcal{H}$, we now conclude that
\[\P(\GC \in \mathcal{A}) \leq \P(B^c) + \P(Q^c) + \P(\TC \in \mathcal{H}) + (\P(\tGC \in \mathcal{A}))^{1/2}\,.\]
Taking a limit, we get that $\lim\sup_{n \to \infty}\P(\GC \in \mathcal{A}) \leq \delta$ and the proof is completed by letting $\delta \to 0$.
%
\end{proof}

\section{The early giant component}\label{sec:simplify-model}

In this section, we consider the special case of Theorem~\ref{mainthm-struct-gen} for $\epsilon=o(n^{-1/4})$, and namely prove
Theorem~\ref{mainthm-struct}. We show how each of the three steps described in Theorem~\ref{mainthm-struct-gen} reduces to the corresponding steps in Theorem~\ref{mainthm-struct}.

\subsection{Step~\ref{item-struct-base}: The kernel}
Let $\Lambda \sim \mathcal{N}\sim(1+\epsilon-\mu, \frac{1}{\epsilon n})$. By \eqref{eq-mu-taylor}, we have  $\mu = 1- \epsilon + O(\epsilon^2)$, and so
\begin{align*}
 \E \Lambda &= 1 + \epsilon - \mu = 2\epsilon + O(\epsilon^2)~,\quad \var (\Lambda) = \frac{1}{\sqrt{\epsilon n}}\,,
\end{align*}
giving that $\Lambda = (2+o(1))\epsilon$ w.h.p. In particular, the probability that $D_u \geq 4$ for some vertex $u$ is
\[ \P(\Po(\Lambda) \geq 4) = O(\epsilon^4) = o(n^{-1})\,, \]
and a union bound thus implies that the kernel is 3-regular w.h.p. In other words, we have $N_k = 0$ for all $k \geq 4$, and
\[ N = N_3 \sim \Bin(n, \mathrm{e}^{-\Lambda} \Lambda^3/6 )\mbox{ conditioned to be even.}\]

It remains to compare the distributions of $N$ and $2\lfloor Z\rfloor$, where
\[Z \sim \mathcal{N}(\tfrac23\epsilon^3 n, \epsilon^3 n)\,.\]
The first step in this direction is to approximate the binomial variable by a Poisson variable. A well-known and straightforward application of the Stein-Chen method (see, e.g., \cite{BC}) is that for any $n$ and $q$,
\[ \| \Bin(n,q) - \Po(nq) \|_\mathrm{TV} \leq  q \;\wedge\; nq^2\,,\]
where the total-variation distance $\|\cdot\|_\mathrm{TV}$ between two distributions $\sigma,\pi$ on a finite space $\Omega$ is given by
\begin{equation}\label{eq-tv-def} \|\sigma-\pi\|_\mathrm{TV} \deq \sup_{A \subset \Omega} \left|\sigma(A)-\pi(A)\right| = \frac12\sum_{x\in\Omega} \left|\sigma(x)-\pi(x)\right|\,.\end{equation}
Therefore, given that $\Lambda = (2+o(1))\epsilon$ (again, this holds w.h.p.) we have
\[ \| \Bin(n,\mathrm{e}^{-\Lambda}\Lambda^3/6) - \Po(n \mathrm{e}^{-\Lambda}\Lambda^3/6) \|_\mathrm{TV} \leq  O(n \epsilon^6) = o(n^{-1/2})\,.\]
Clearly, by definition \eqref{eq-tv-def}, a negligible total-variation distance between two distributions already implies they are contiguous (in both directions), hence it suffices to compare $2\lfloor Z\rfloor$ to the variable $Y$, distributed as $\Po(n \mathrm{e}^{-\Lambda}\Lambda^3/6)$ conditioned to be even.
We will show that for some region $\mathcal{Q}$ such that $\P(Y\in \mathcal{Q}) \to 1$ and some $c = c(\mathcal{Q}) >0$,
\begin{equation}\label{eq-compare-Y-Z} \P(Y = t) \leq c \cdot \P(2\lfloor Z\rfloor = t)\quad \mbox{ for all even $t\in\mathcal{Q}$}\,.\end{equation}
Let $\delta > 0$. By the above properties of $\Lambda$, there exists some $M > 0$ such that
\[ \P\left(|\Lambda - 2\epsilon| > M/\sqrt{\epsilon n}\right) < \delta\,.\]
The following result is a special case of a theorem of \cite{Durrett}.
\begin{theorem}[\cite{Durrett}*{Ch.\ 2, Theorem 5.2}, reformulated]
Let $X$ be a random variable on $\N$ with $\P(X = k) > 0$ for all $k\in\N$. Suppose that $\E X = \nu< \infty$ and $\var X = \sigma^2 < \infty$.
Let $X_i$ be i.i.d.\ distributed as $X$ and $S_m = \sum_{i=1}^m X_i$. Then as $m \to \infty$, we have
$$\sup_{x \in \mathcal{L}_m}\left|\sqrt{m}\,\P \left(\frac{S_m -  m\nu}{\sqrt{m}} = x\right) - \frac{1}{\sqrt{2\pi}\sigma} \mathrm{e}^{-x^2/\sigma^2}\right| \to 0~,$$
where $\mathcal{L}_m = \{(z - m \nu )/ \sqrt{m}: z\in \Z\}$.
\end{theorem}
In our setting, given that $|\Lambda-2\epsilon| \leq M/\sqrt{\epsilon n}$ we have a Poisson distribution with parameter
\[m=\tfrac43\epsilon^3 n +  O(\sqrt{\epsilon^3 n})\,,\] and clearly the effect of conditioning that it is even, as well as rounding $m$ to its integer part, only affect the density function by a constant factor. This translates into a sum of $\lfloor m\rfloor$ i.i.d.\ $\Po(1)$ variables, and by the above theorem we conclude that, as long as $|\Lambda-2\epsilon|< M/\sqrt{\epsilon n}$, there exists some $c =c(M)> 0$ such that
\[ \P(Y = t) \leq \frac{c}{\sqrt{\epsilon^3 n}} \quad\mbox{ for any integer $t$}\,. \]
Furthermore, we can choose $\mathcal{Q} = [ m - M'\sqrt{\epsilon^3 n}, m+ M'\sqrt{\epsilon^3 n} ]$ for a suitably large $M'>0$ so that
\[ \P(Y \notin \mathcal{Q}) < \delta \,.\]
By the definition of $Z$ (note that $2\lfloor Z\rfloor$ has mean $\frac43 \epsilon^3 n + O(1)$ and variance of order $\epsilon^3 n$), there exists some $c'=c'(M')$ such that
\[ \P(2\lfloor Z\rfloor = t) \geq \frac{c'}{\sqrt{\epsilon^3 n}}\quad\mbox{ for all even $t \in \mathcal{Q}$}\,.\]
Altogether, it follows that for any sequence of subsets of integers $\mathcal{S}=\mathcal{S}(n)$, if $\P(2\lfloor Z \rfloor \in \mathcal{S}) = o(1)$ then $\P(N \in \mathcal{S}) \leq 2\delta + o(1)$. We now let $\delta \to 0$ to complete the contiguity of $N$ and $2\lfloor Z\rfloor$.

\subsection{Step~\ref{item-struct-edges}: The 2-core}

Here we need to compare the effect of replacing the 2-paths by i.i.d.\ $\Geom(\epsilon)$ variables rather than $\Geom(1-\mu)$.
Rather than just showing contiguity between the two models, we will show a stronger statement, namely that the total-variation between the joint distributions of the path lengths are negligible.
The total-variation distance $\|\cdot\|_\mathrm{TV}$ between two distributions $\sigma,\pi$ on a finite space $\Omega$ is given by
\[ \|\sigma-\pi\|_\mathrm{TV} \deq \sup_{A \subset \Omega} \left|\sigma(A)-\pi(A)\right| = \frac12\sum_{x\in\Omega} \left|\sigma(x)-\pi(x)\right|\,.\]
In our case, since the ratio of $\frac{p(1-p)^k}{q(1-q)^k}$ for $0<p<q<1$ is monotone increasing in $k$, we can clearly consider sets of the $A = \{k,k+1,\ldots\}$ in the supremum above. Hence,
\begin{align*}
\|\Geom(1-\mu) &- \Geom(\epsilon)\|_{\mathrm{TV}} = \sup_k |\mu^k - (1 - \epsilon)^k|\\
 &\leq \sup_k |\mu - (1 -\epsilon)| \cdot k [\mu \vee (1-\epsilon)]^{k-1}   = O(\epsilon)\,,
\end{align*}
where, as $|\mu-(1-\epsilon)| = \Theta(\epsilon^2)$, the value of $k$ optimizing the above has order $1/\epsilon$.

Recalling that $|E(\K)| = O(\epsilon^3 n)$ w.h.p., we infer that the total-variation distance between the joint distribution
of the 2-path lengths in the two models (i.i.d.\ $\Geom(1-\mu)$ variables and i.i.d.\ $\Geom(\epsilon)$ variables) is
$O(\epsilon^4 n)$. Our assumption that $\epsilon=o(n^{-1/4})$ now gives that this is $o(1)$.

\subsection{Step~\ref{item-struct-bushes}: The attached trees}
We now wish to compare the distributions of i.i.d.\ PGW($\mu$)-trees to i.i.d.\ PGW($1-\epsilon$)-trees.
Recall that the size of a PGW($\gamma$)-tree follows a Borel($\gamma$) distribution, as given in \eqref{eq-PGW-size}. Thus,
\begin{align*}
\|\mathrm{Borel}(\mu) - \mathrm{Borel}(1 - \epsilon)\|_{\mathrm{TV}}& = \frac{1}{2}\sum_{t} \frac{t^{t-1}}{t!} \big|\mu^{t-1} \mathrm{e}^{-\mu t} - (1 - \epsilon)^{t - 1} \mathrm{e}^{-(1 - \epsilon)t}\big|\\
&\leq O(\epsilon^2) \sum_t \frac{t^{t-1}}{t!} |x_t^{t-2} \mathrm{e}^{-x_t t}(t - 1 -tx_t)|\,,
\end{align*}
 Noticing that $x_t = 1 - \epsilon +O(\epsilon^2)$ and applying Stirling's inequality, we get
\begin{align*}\frac{t^{t-1}}{t!}&\left|x_t^{t-2}\mathrm{e}^{-x_t t}(t-1-tx_t) \right| = \Theta(1)
t^{-3/2} (x_t \mathrm{e}^{1-x_t})^t \left|-1+t\epsilon + O(t\epsilon^2)\right| \\
& \leq O(t^{-3/2}) \mathrm{e}^{-(1-x_t)^2/2} (1+ t\epsilon)
\leq O(t^{-3/2}) \mathrm{e}^{-\epsilon^2t/3} (1+ t\epsilon) \,,
\end{align*}
where the first inequality is by the fact that $1 - y \leq \mathrm{e}^{-y - y^2/2}$ for all $y>0$, and the second one holds for any large $n$ by the definition of $x_t$. Therefore,
\begin{align*}\|\mathrm{Borel}(\mu) &- \mathrm{Borel}(1 - \epsilon)\|_{\mathrm{TV}} = O(\epsilon^2) + O\bigg(\epsilon^2 \sum_t  \frac{\epsilon \mathrm{e}^{-\epsilon^2t/3}}{\sqrt{t}}\bigg) \\
& \leq O(\epsilon^2) + O(\epsilon^2) \int_0^\infty \frac1{\sqrt{\epsilon^2 x/3}} \mathrm{e}^{-\epsilon^2 x/3 } d(\epsilon^2x/3) =O(\epsilon^2)\,,
\end{align*}
where we used the fact that $\int_0^\infty \frac1{\sqrt{y}}\mathrm{e}^{-y}$ converges (to $\sqrt{\pi}$).

With high probability, the size of the 2-core (that is, the number of PGW-trees we attach) is $O(\epsilon^2 n)$, and so the total-variation distance between the two joint distributions is at most  $O(\epsilon^4 n) = o(1)$, as required.

\section{Analysis of the Cut-Off Line Algorithm}\label{sec:lambda-C}
In this section, we analyze Algorithm~\ref{algorithm-cola} for generating the Poisson cloning model, and establish a tight concentration result for $\Lambda_C$ (the location of the cut-off line when all light clones are exhausted), as stated in Theorem~\ref{thm-Lambda-C-upper}.

\begin{proof}[\textbf{\emph{Proof of Theorem~\ref{thm-Lambda-C-upper}}}]
We wish to prove inequality~\eqref{eq-Lambda-C-window}, i.e., that for some fixed $c>0$, the probability that
$|\Lambda_C - \theta_\lambda \lambda| \geq \tfrac{\gamma}{\sqrt{\theta_\lambda n}}$ is at most $\exp(-c \gamma^2)$.

Notice that, prior to the first time the algorithm reaches Step~\ref{item-cola-2}, the notion of active/passive vertices does not play a role in its decisions. Since this is the only change between subsequent phases, it follows that $\Lambda_C$ is precisely the same regardless of the choice of phase boundaries.
In particular, we may choose $\beta$ as follows: Take $\frac{1 - \theta_\lambda}{3} \leq \beta \leq \frac{1-\theta_\lambda}{2}$ and an integer $m$ such that
\begin{equation}
  \label{eq-beta-m-choice}
  (1-\beta)^{m-1} =
\theta_{\lambda} + \frac{\gamma}{\sqrt{\theta_{\lambda} n}}~,
\end{equation} where $\gamma$ as given in the lemma, i.e., $\gamma = o\big(\sqrt{\theta_\lambda^3 n}\big)$.

In order to prove the lemma, we first estimate the number of $j$-active clones for each $j$, denoted by $N_j$.
Let $M_j$ be the number of $j$-active clones that are \emph{matched} during phase $j$. We need the following lemma to estimate $M_j$ given $N_j$.

\begin{lemma}[\cite{KimA}*{Lemma 2.2}]\label{lem-cut-off-line}
 Consider a Poisson $\mu$-cell for $\mu > 0$, and let $N$ be its total number of clones. For $0 < \theta \leq 1$, let $N(\theta)$ be the number of matched clones once the cut-off line reaches $\theta \mu$. Then there exists some $c > 0$ so that the following holds: For any $0 <\theta_0 < 1$  and $l,\Delta>0$,
$$\P\left(\max_{\theta_0 \leq \theta \leq 1} |N(\theta) - (1 - \theta^2)k| \geq \Delta \,\Big|\, N = k\right) \leq 2\exp\big[-c\big(\Delta \wedge \tfrac{\Delta^2}{(1-\theta_0)k}\big)\big]~.$$
\end{lemma}

By definition of the Cut-Off Line Algorithm, if either one of the two
unmatched clones of a passive vertex was matched in a given step, then the other clone is guaranteed to be matched in the next step (either in this phase or in a later one), as it is inserted to the top of the stack. This means that, for the purpose of determining the number of matched active clones throughout the phase, $M_j$, applying the algorithm with or without passive vertices is effectively the same (since one can always identify the two clones of each passive vertex, then contract the 2-paths into edges between active clones).

That said, one must consider the following delicate point. If the end-of-phase boundary is reached while the top of the stack contains
a passive clone, this corresponds to a path whose one endpoint is an active clone $(u,i)$, yet its other endpoint is a passive clone.
In this case, the active clone should not be considered as matched when disregarding all passive clones. Let $A_j$ denote this event for phase $j$, and define $$M'_j \deq M_j - \one_{A_j}$$ to be the number of $j$-active clones that are matched during phase $j$ while disregarding (contracting) the passive clones.

Combining the above observation with the fact that phase $j$ began
with $N_j$ active clones and a cut-off line at $\mu=(1-\beta)^{j-1}\lambda$ and ended as soon as the cut-off line reached $(1-\beta)\mu$,
we apply Lemma \ref{lem-cut-off-line} for $\theta = \theta_0 = 1-\beta$ and conclude that for some constant $c>0$,
\begin{equation}\label{eq-bound-M-j}
\P\left(\left|M'_j - (1- (1- \beta)^2)N_j\right| \geq \Delta \given N_j\right) \leq 2\exp\big[-c\big(\Delta \wedge \tfrac{\Delta^2}{\beta N_j}\big)\big]~.
\end{equation}
Let $B_j$ denote the number of vertices which have precisely $2$ clones to the left of the end-of-phase boundary of phase $j$, and at least $1$ more clone in the interval of phase $j$. Note that, for such a vertex $v$, it is clearly $j$-active, and it would become $(j+1)$-passive if and only if the formerly mentioned $2$ clones are unmatched by the end of phase $j$. In this case, two formerly active clones will be relabeled as passive. In particular, the number of clones that transition from being $j$-active to being $(j+1)$-passive is at most $2B_j$.

On the other hand, a clone can transition from being $j$-passive to being $(j+1)$-active if and only if it happened to be at the top of the stack when phase $j$ ended, and in particular, the event $A_j$ occurred.

Adding these two, along with the number of $j$-active clones matched in this phase $M'_j$, we conclude that the number of $(j+1)$-active clones satisfies
\begin{equation}\label{eq-N-M-B}
N_{j+1} \geq N_j - M_j + \one_{A_j} - 2 B_j = N_j - M'_j - 2 B_j~.
\end{equation}
Note that, as long as Step~\ref{item-cola-2} of the algorithm is not reached, the stack always consists of light clones exclusively. Therefore, up till that point, if a vertex has $2$ unmatched clones, both will remain unmatched until the cut-off line reaches one of them. Suppose that phase $j_0$ is the first one where the algorithm invoked Step~\ref{item-cola-2}. In that case, for any $j < j_0$, the vertices counted in $B_j$ are precisely those that were $j$-active yet became $(j+1)$-passive. We deduce that \eqref{eq-N-M-B} is in fact an equality for all $j < j_0$.

In order to analyze $B_j$, for each $v \in V$ and $0 \leq \theta < \theta' \leq 1$ let $d_v(\theta, \theta')$ denote the
number of $v$-clones whose assigned value belongs to $[\theta\lambda, \theta'\lambda)$. Further let  $d_v(\theta) \deq d_v(0, \theta)$. Recall that phase $j$ begins with the line at $(1-\beta)^{j-1}\lambda$ and ends with the line at $(1-\beta)^j\lambda$.
Hence, for $\theta_j = (1-\beta)^{j-1}$, we have by the definition of $B_j$ that
$$B_j = \sum_{v\in V} \one_{\{d_v(\theta_{j+1}) = 2\}}\ \one_{\{d_v(\theta_{j+1}, \theta_j) \geq 1\}}~.$$
Observe that $(d_v(\theta_{j+1}), d_v(\theta_{j+1}, \theta_j))$ for $v\in V$ are i.i.d.\ pairs of independent Poisson random variables
with means $\theta_{j+1} \lambda$ and $(\theta_j - \theta_{j+1}) \lambda$ respectively. Applying Chernoff's bound (cf., e.g., \cite{AS}) we have that for some $c_1 > 0$,
\begin{equation}\label{eq-bound-B-j}
\P\left(\Big| B_j - \frac{(\theta_{j+1}\lambda)^2}{2} \mathrm{e}^{-\theta_{j+1}\lambda}(1-\mathrm{e}^{-\beta \theta_j \lambda}) n\Big| \geq \Delta\right)
\leq 2 \exp\Big(-c_1 \tfrac{\Delta^2}{\theta_j^3 n}\Big)~.
\end{equation}
Combined with \eqref{eq-bound-M-j} and \eqref{eq-N-M-B}, we arrive at the following estimated lower bound for $N_{j+1}$:
$$(1-\beta)^2 N_j - (\theta_{j+1}\lambda)^2
\mathrm{e}^{-\theta_{j+1}\lambda} (1 - \mathrm{e}^{-\beta \theta_j \lambda}) n~.$$
Applying this inductively, we expect that the following would be a lower bound for $N_j$:
\begin{equation}\label{eq-nj-guess}
\theta_j^2 \lambda (1 - \lambda \mathrm{e}^{-\theta_j \lambda}) n
\end{equation}
(indeed, this is later shown in Lemma~\ref{lem-concentration-N-L-lower}).

Suppose that the algorithm is at the beginning of phase $j$, and Step~\ref{item-cola-2} has not been reached yet (in any of the phases thus far). By the discussion above, any $j$-active vertex has either $1$ or strictly more than $2$ clones with values in $(0,\theta_j\lambda)$. In particular, the number of $j$-active clones that are \emph{heavy} at the beginning of phase $j$ can then be written as
$$H_j = \sum_{v\in V} d_v(\theta_j) \one_{\{d_v(\theta_j) > 2\}}~,$$
and the number of light clones at the start of phase $j$ is then precisely
$$L_j \deq N_j - H_j~.$$
In general (once Step~\ref{item-cola-2} is invoked), $H_j$ is an upper bound for the number of $j$-active clones that are heavy at this point. (Note that the only reason for this bound not to be tight is on account of clones that are already matched. That is, $\sum_{v\in V} \one_{\{d_v(\theta_j) > 2\}}$ counts all $j$-active heavy \emph{vertices}, in addition to perhaps some whose clones are all matched by phase $j$.)
Hence, $L_j$ is always a lower bound for the number of light
clones at the start of phase $j$.

We need the following large deviation inequality:
\begin{lemma}[\cite{KimB}*{Corollary 4.2}]\label{lem-Generalized-Chernoff-bound}
Let $X_1, \ldots, X_m$ be a sequence of independent random variables. Suppose $\E[X_i] = \mu_i$ and there are $b_i$, $d_i$ and $\xi_0$ such that $\E[(X_i - \mu_i)^2] \leq b_i$, and
$$\left|\E\left[(X_i - \mu_i)^3 \mathrm{e}^{\xi(X_i - \mu_i)}\right]\right| \leq d_i \qquad \mbox{for all } 0\leq |\xi| \leq \xi_0~.$$
If $\delta \xi_0 \sum_{i=1}^m d_i \leq \sum_{i=1}^{m} b_i$ for some $0 < \delta \leq 1$, then
$$\P \left(\big|\mbox{$\sum_{i=1}^m X_i - \sum_{i=1}^m \mu_i$}\big| \geq \Delta\right) \leq \mathrm{e}^{-\frac{1}{3}\min\{\delta \xi_0 \Delta, \frac{\Delta^2}{\sum_{i=1}^m b_i}\}}~,$$
for all $\Delta>0$.
\end{lemma}
Since $d_v(\theta_j)$ are i.i.d.\
$\Po(\theta_j \lambda)$ variables, an application of the above lemma gives that for some constant $c_2>0$,
\begin{equation}\label{eq-bound-H-j}
\P\left(|H_j - \theta_j \lambda(1 - \mathrm{e}^{-\theta_j \lambda} - \theta_j \lambda \mathrm{e}^{-\theta_j \lambda}) n| \geq \Delta\right) \leq
2\exp\big[-c_2\big(\Delta \wedge \tfrac{\Delta^2}{\theta_j^3 n}\big)\big]~.
\end{equation}
Recalling that $m$ is an integer with $(1-\beta)^{m-1} =
\theta_{\lambda} + \frac{\gamma}{\sqrt{\theta_{\lambda} n}}$ (see \eqref{eq-beta-m-choice}), set
\begin{equation}\label{eq-Delta-j-def}
\Delta_j \deq \frac{\gamma}{100} \sqrt{\theta_j^3 n} \sum_{i=1}^j (1-\beta)^{(2j-i-m)/4}~.
\end{equation}
Observe that the following sequence is increasing in $j$:
$$(1-\beta)^{\frac{j}{4}} \theta_j^{-3/2} = (1- \beta)^{3/2} (1-\beta)^{-5j/4}~.$$
We then get that for all $j \in [m]$,
\begin{equation}\label{eq-bound-gamma-theta}
\gamma (1 - \beta)^{(j-m)/4} (\theta_j^3  n)^{-1/2} \leq \frac{\gamma}{\sqrt{\theta_m^3 n}} = o(1)~,
\end{equation}
where the last equality used the facts $\gamma = o\big(\sqrt{\theta_\lambda^3 n}\big)$ and
\begin{equation}\label{eq-theta-m-asymp}\theta_m = (1-\beta)^{m-1} = \theta_\lambda + \frac{\gamma}{\sqrt{\theta_\lambda n}} = (1 + o(1))\theta_\lambda~.\end{equation}
It then follows from \eqref{eq-bound-gamma-theta} that for any $j\in [m]$,
\begin{align}
\Delta_j& = \frac{\gamma}{100} \sqrt{\theta_j^3 n} \sum_{i=1}^j (1-\beta)^{(2j-i-m)/4}\nonumber\\
 &= \frac{\gamma}{100} \sqrt{\theta_j^3 n} (1-\beta)^{(j-m)/4} \sum_{i=1}^j (1-\beta)^{(j-i)/4} = o (\theta_j^3 n)~.\label{eq-Delta-j-bound}
\end{align}
Recalling \eqref{eq-nj-guess} and \eqref{eq-bound-H-j}, we will next establish lower bounds for the $N_j$'s and $L_j$'s in terms of the following parameters:
\begin{align}
\begin{array}
  {rcl}n_j &\deq& \theta_j^2 \lambda (1 - \lambda \mathrm{e}^{-\theta_j \lambda})n\\
  l_j &\deq& \theta_j \lambda (\theta_j -1 + \mathrm{e}^{-\theta_j \lambda}) n
\end{array} &\qquad\mbox{for $j\in[m]$}\,.\label{eq-nj-lj-def}
\end{align}

\begin{lemma}\label{lem-concentration-N-L-lower}
There exists a constant $c>0$ such that the following holds:
\begin{align*}
\P\left(\exists j \in [m]: N_j < n_j - \Delta_j\right) \leq \mathrm{e}^{-c \gamma^2 }~,\\
\P \left(\exists j\in [m]: L_j < l_j - 2 \Delta_j\right) \leq \mathrm{e}^{-c\gamma^2}~.
\end{align*}
\end{lemma}
\begin{proof}
With \eqref{eq-bound-B-j} in mind, define the following for each $j \in  [m]$:
\begin{align}
 b_j &\deq \frac{(\theta_{j+1}\lambda)^2}{2}
\mathrm{e}^{-\theta_{j+1}\lambda}(1-\mathrm{e}^{-\beta \theta_j \lambda})n~,\nonumber\\
\gamma_j& \deq (1-\beta)^{(j-m)/4} \gamma~.\label{eq-def-gamma-j}
\end{align}
It is clear (see definition \eqref{eq-nj-lj-def}) that
\begin{equation}\label{eq-n-b-j}
n_{j+1} = (1-\beta)^2 n_j - 2 b_j~,
\end{equation}
and furthermore, by \eqref{eq-bound-gamma-theta}, we have that
\begin{equation}\label{eq-gamma-j}
\gamma_j = o \big(\sqrt{\theta_j^3 n}\big)~.
\end{equation}
For $\ell \in [m]$, decomposing the events in the required lower bound on $N_j$ gives
\begin{align*}
\P\left(\exists j\in [\ell]: N_j < n_j - \Delta_j\right)=\P\big(\exists j \in [\ell-1]: N_j < n_j - \Delta_j&\big)\\
+~ \P\big(\cap_{j\in[\ell-1]} \big\{N_j \geq n_j - \Delta_j\big\} \cap \big\{N_{\ell} < n_{\ell} - \Delta_{\ell}\big\}&\big)~,
\end{align*}
as well as
\begin{align*}
\P\big(\cap_{j\in[\ell-1]} &\big\{N_j \geq n_j - \Delta_j\big\} \cap \big\{N_{\ell} < n_{\ell} - \Delta_{\ell}\big\}\big)\\
&\leq \P \left(N_{\ell-1} \geq n_{\ell-1}- \Delta_{\ell-1}~,~ N_{\ell} < n_{\ell} - \Delta_{\ell}\right) \deq P_{\ell}~.
\end{align*}
Recall that $\gamma_j$ is decreasing in $j$, thereby for any constant $c_1'>0$ there exists a constant $c_2'>0$ such that
$$\sum_{j=1}^{\ell} \mathrm{e}^{-c_1'\gamma_j^2} \leq
\mathrm{e}^{-c_2' \gamma_{_{\ell}}^2}~.$$
It will thus suffice to show that $P_{j} \leq \mathrm{e}^{-c'
\gamma_{j}^2}$ for some constant $c' > 0$ and every $j\in[m]$.

For $j=1$, recall that $N_1$ is the number of active clones at the beginning of the algorithm
(since all clones are initially unmatched, the passive vertices are those with precisely $2$ clones, and all other vertices are active)
and $n_1 = \lambda(1-\lambda\mathrm{e}^{-\lambda})n = \E N_1$. In addition, $\Delta_1 = \frac{1}{100}\gamma_1 \sqrt{n}$
and $\gamma_1\to\infty$ with $n$ (by the fact that $(1-\beta)^{-m} \asymp 1/\epsilon \to\infty$). Hence, by a standard application of the Central Limit Theorem to the i.i.d.\ random variables defined by the number of active clones that each vertex contributes, we deduce that $P_1 \leq \mathrm{e}^{-c_0\gamma_1^2}$ for some $c_0>0$ fixed.

Next, consider $P_{j+1}$ for $j\in\{1,\ldots,m-1\}$. Combining \eqref{eq-N-M-B} with \eqref{eq-n-b-j} we get that for each such $j$
\begin{align*}
N_{j+1} -n_{j+1} &\geq N_j - (1-\beta)^2 n_j -  M'_j - 2(B_j - b_j) \\
&= (1-(1-\beta)^2)N_j - M'_j + (1-\beta)^2(N_j - n_j) - 2(B_j - b_j)~.
\end{align*}
Therefore, the event addressed in $P_{j+1}$ implies that
$$\Delta_{j+1} < n_{j+1} - N_{j+1} \leq M'_j - (1-(1-\beta)^2)N_j + (1-\beta)^2 \Delta_j + 2(B_j - b_j)~,$$
and also $N_j \geq n_j - \Delta_j$. In particular, $P_{j+1}$ is at most the probability that
\begin{align*}
&M'_j - (1-(1-\beta)^2)N_j + 2(B_j-b_j) \geq \Delta_{j+1} - (1-\beta)^2\Delta_j~,\mbox{ and }\\
&N_j \geq n_j - \Delta_j~.
\end{align*}
Adding the fact that, by definitions \eqref{eq-Delta-j-def} and \eqref{eq-def-gamma-j},
 $$\Delta_{j + 1} = (1-\beta)^2 \Delta_j + \frac{\gamma_{j + 1}}{100}
\sqrt{\theta_{j + 1}^3 n}~,$$
and we deduce that
\begin{align*}
P_{j+1} 
&\leq \P\left(B_j - b_j > \frac{\gamma_{_{j+1}}}{400}\sqrt{\theta_{j+1}^3 n}\right)\\
& + \P\left(M'_j - (1-(1-\beta)^2)N_j > \frac{\gamma_{_{j+1}}}{200} \sqrt{\theta_{\ell + 1}^3 n}\given N_j \geq n_j - \Delta_j \right)~.
\end{align*}
Combining \eqref{eq-bound-B-j} and \eqref{eq-gamma-j}, we can obtain an upper bound $\mathrm{e}^{-c_1 \gamma_{\ell+1}^2}$ on the first term, for
some constant $c_1 > 0$. At the same time, \eqref{eq-bound-M-j} provides an upper bound of $\mathrm{e}^{-c_2 \gamma_{_{\ell+1}}^2}$ on the second term, for some $c_2 > 0$ fixed.

Altogether, we have shown the desired upper bound on $P_j$ for all $j\in[m]$, implying the first inequality in the lemma.

Recall now that inequality \eqref{eq-bound-H-j} gives that for some constant $c_3,c_4 > 0$,
\begin{align*}
&\P\left(\exists j\in [\ell]: |H_j - \theta_j \lambda(1 - \mathrm{e}^{-\theta_j \lambda} - \theta_j \lambda
\mathrm{e}^{-\theta_j \lambda}) n| \geq \Delta_j\right)\\
&\leq 2\sum_{j=1}^{\ell} \mathrm{e}^{-c_3 \gamma_j^2} \leq \mathrm{e}^{-c_4\gamma_\ell^2}~.
\end{align*}
Combining this with the fact $L_j \geq N_j - H_j$, as well as the above lower bound on $N_j$, yields the second statement of the lemma, as required.
\end{proof}

We can now derive a lower bound on the first time that Step~\ref{item-cola-2} is applied (that is, the first
time at which there are no light clones). Equivalently, this gives an upper bound on $\Lambda_C$ (the $x$-coordinate
of the cut-off line at that point).

By the definition, the number of light clones throughout the algorithm has the following property:
\begin{itemize}
\item As long as there light clones in the stack, in
each step one of them will be popped and matched, and as a result, at most one new light clone will be created.
\item  If there are no light clones in the stack (and the algorithm is not concluded) then following Step~\ref{item-cola-2} the stack
will necessarily be comprised of a single heavy clone. This clone will then be popped in the next iteration of Step~\ref{item-cola-1},
while creating at most two new light clones.
\end{itemize}
That is to say, once the number of light clones drops to $0$, it can never again exceed $2$.
In particular, if all the light clones disappear for the first time during phase $j$ for some $j = 1,\ldots, m-1$, we must have that $L_m \leq 2$ (since $L_j$ is a lower bound on the number of light clones at the start of phase $j$).

By \eqref{eq-nj-lj-def} we have that
\begin{align*}
l_m &= \theta_m \lambda(\theta_m -1 + \mathrm{e}^{-\theta_m \lambda} )n~.
\end{align*}
By the definition of $\theta_\lambda$ and its asymptotic behavior (see \eqref{eq-theta-lambda-def},\eqref{eq-theta-lambda}), as well as the fact
that $\theta_m=(1+o(1))\theta_\lambda$, we have that for all $\theta_\lambda \leq x \leq \theta_m$ the function $f(x) = x - 1 + \mathrm{e}^{-\lambda x}$ satisfies
$$ f'(x) = 1 - \lambda \mathrm{e}^{-\lambda x} = 1 - \lambda \mathrm{e}^{- (1+o(1))\theta_\lambda \lambda} =
1-\lambda(1-\theta_\lambda)(1-o(\theta_\lambda)) = (1-o(1))\epsilon~.$$
Since $f(\theta_\lambda)=0$ and $\theta_m = \theta_\lambda + \frac{\gamma}{\sqrt{\theta_m n}} = (2+o(1))\epsilon$ (see \eqref{eq-theta-m-asymp}), we can apply the Mean Value Theorem and get
\begin{equation}
  \label{eq-lm-order}
 l_m = \theta_m \lambda f(\theta_m)  n = \theta_m \lambda \cdot (1-o(1))\epsilon \frac{\gamma}{\sqrt{\theta_m n}} n
= (\tfrac12-o(1))\gamma\sqrt{\theta_m^3 n}~.
\end{equation}
On the other hand, by \eqref{eq-Delta-j-def} (and recalling requirement \eqref{eq-beta-req1} from $\beta$)
\begin{align}
 \Delta_m  &= \frac{\gamma}{100}\sqrt{\theta_m^3 n} \sum_{i=1}^m (1-\beta)^{(m-i)/4}
  \leq \frac{\gamma}{100}\sqrt{\theta_m^3 n} \frac{1}{1-(1-\beta)^{1/4}} \nonumber\\
  &\leq \frac{\gamma}{100}\sqrt{\theta_m^3 n} \frac{1}{1-(\frac23 - o(1))^{1/4}} \leq \frac{\gamma}8 \sqrt{\theta_m^3 n}
  \leq \tfrac13 l_m~,  \label{eq-Deltam-order}
  \end{align}
where the last two inequalities hold for any large $n$.
As $\theta_m = (1+o(1))\theta_\lambda$ and $\theta_\lambda^3 n \to\infty$ with $n$, we immediately have
that $l_m\to\infty$ as well. However, Lemma~\ref{lem-concentration-N-L-lower} gives that
$$ \P(L_m \leq 2) \leq \P(L_m \leq l_m - 2 \Delta_m) \leq \mathrm{e}^{-c\gamma^2} $$
for some fixed $c > 0$. By the above discussion, this translates into an upper bound of $\Lambda_C$:
\begin{equation}
\P\left(\Lambda_C \geq \theta_\lambda \lambda + \tfrac{\gamma}{\sqrt{\theta_\lambda n}}\right) \leq
\mathrm{e}^{-c\gamma^2}~.\label{eq-Lambda-C-upper-bound}
\end{equation}

\bigskip

The above upper bound on $\Lambda_C$ ensures that
Step~\ref{item-cola-2} is not applied in the first $m-1$ phases, except with probability $\exp(-c\gamma^2)$.
Recall that, if Step~\ref{item-cola-2} has not yet been applied in phases $1,\ldots,j-1$ then our lower bound \eqref{eq-N-M-B} for $N_j$ is in fact
an equality, and similarly, $L_j$ is precisely the number of light clones at the beginning of phase $j$.
Therefore, assuming that indeed Step~\ref{item-cola-2} was not applied in phases $1,\ldots,m-1$ (we account for the above error probability), we may now apply the same proof of Lemma~\ref{lem-concentration-N-L-lower}, this time
with respect to the events $\{N_j > n_j + \Delta_j\}$ and $\{L_j > l_j + 2 \Delta_j\}$.
This gives the following matching upper bounds on the $N_j$'s and $L_j$'s:
\begin{lemma}\label{lem-concentration-N-L-upper}
There exists a constant $c>0$ such that the following holds:
\begin{align*}
\P\left(\exists j \in [m]: N_j > n_j + \Delta_j\right) \leq \mathrm{e}^{-c \gamma^2 }~,\\
\P \left(\exists j\in [m]: L_j > l_j + 2 \Delta_j\right) \leq \mathrm{e}^{-c\gamma^2}~.
\end{align*}
\end{lemma}

To obtain the required upper bound on $\Lambda_C$, assume that Step~\ref{item-cola-2} was not applied
in phases $1,\ldots,m-1$ (this happens except with probability $\exp(-c\gamma^2)$). We now wish to show that
all light clones will disappear shortly after commencing phase $m$ with probability at least $1-\exp(-c\gamma^2)$.

Since we did not apply Step~\ref{item-cola-2} yet, the stack exclusively contains light clones,
and a clone is active (more precisely, $m$-active) iff it has $3$ unmatched clones or more.
Now, if we ignore the passive clones, the algorithm must remove at least $2$ light clones from the stack in order to create
a new light clone.

Suppose that at the beginning of phase $m$, the stack contains $k$ light clones. By the above discussion, after matching $k$ light
clones, the stack will be of size at most $k/2$. Iterating, it follows that
after matching at most $2k$ active clones, every light clone will disappear (the stack will be exhausted).

By Lemma~\ref{lem-concentration-N-L-upper} and the fact that $l_m+2\Delta_m \leq 2 \gamma\sqrt{\theta_m^3 n}$
for any large $n$ with room to spare (as established in \eqref{eq-lm-order},\eqref{eq-Deltam-order}), there are
$L_m \leq 2\gamma\sqrt{\theta_\lambda^3 n}$ light clones
at the beginning of phase $m$, except with probability $1- \mathrm{e}^{-c\gamma^2}$.
Combined with the above argument, we conclude that the stack of light clones will be exhausted after
 matching at most $4\gamma\sqrt{\theta_\lambda^3 n}$ active clones, except with the above error probability.

We will next show that at least this many active clones will be matched by the time
the cut-off line reaches the point $\theta_m \lambda - 10\frac{\gamma}{\sqrt{\theta_\lambda n}}$.
Since $\theta_m \geq \theta_\lambda$, by definition~\eqref{eq-nj-lj-def} we have
\begin{align*}
 n_m &= \theta_m^2 \lambda (1-\lambda\mathrm{e}^{-\theta_m \lambda})n \geq \theta_m^2 \lambda (1-\lambda\mathrm{e}^{-\theta_\lambda \lambda})n
 = \theta_m^2 (1-\lambda(1-\theta_\lambda))n \\
 &= \theta_m^2 ( -\epsilon +(2+o(1))\epsilon) n = (\tfrac12+o(1))\theta_m^3 n~.
\end{align*}
Together with Lemma~\ref{lem-concentration-N-L-lower}, we deduce that for a sufficiently large $n$,
there are $N_m \geq \frac{1}{3}\theta_\lambda^3 n$ unmatched active clones
at the beginning of phase $m$, except with probability $1- \mathrm{e}^{-c\gamma^2}$.
Note that, by the assumption on $\gamma$,
$$\theta_m \lambda - 10\frac{\gamma}{\sqrt{\theta_\lambda n}} = \Big(1 - 10\frac{\gamma}{(1+o(1))\sqrt{\theta_\lambda^3 n}}\Big)\theta_m = (1-o(1))\theta_m~.$$
Since the boundary marking the end of phase $m$ is at $(1-\beta)\theta_m$ (and $\beta$ is bounded away from $0$),
the cut-off line moves through the entire interval between $\theta_m\lambda$ and the above point as part of phase $m$.
Hence, we can use the original version of the Cut-Off Line Algorithm in order to analyze the number of active clones that are matched
along this interval (without considering a potential change of phase).

Applying Lemma~\ref{lem-cut-off-line} with $\theta = 1 - 10\frac{\gamma}{\sqrt{\theta_\lambda^3 n}}$ and $k = \frac{1}{3}\theta_\lambda^3 n$, we can now deduce that, except with probability $\exp(-c\gamma^2)$ we match at least
$$ \tfrac56 (1-\theta^2)k = \tfrac56 (20-o(1))\frac{\gamma}{\sqrt{\theta_\lambda^3 n}}k \geq (5-o(1))\gamma\sqrt{\theta_\lambda^3 n}$$
active clones. Therefore,
\begin{equation*}
\P\left(\Lambda_C \leq \theta_\lambda \lambda - \tfrac{\gamma}{\sqrt{\theta_\lambda n}}\right) \leq
\mathrm{e}^{-c\gamma^{2}}~.
\end{equation*}
Combining this bound with \eqref{eq-Lambda-C-upper-bound} completes the proof of \eqref{eq-Lambda-C-window}.
\end{proof}

\section*{Acknowledgments}
We wish to thank Asaf Nachmias for helpful discussions at an early stage of this project, as well as Nick Wormald for useful comments.

\end{document}